\documentclass[12pt]{amsart}
\usepackage{graphicx,amsmath,amsfonts,latexsym,amssymb,amsthm}
\usepackage[arrow,matrix,curve]{xy}
\usepackage{pictexwd}
\usepackage{dcpic}
\usepackage{verbatim}
\usepackage{a4wide}

\newcommand{\field}[1]{\mathbb{#1}}
\newcommand{\rz}{\field{R}}
\newcommand{\cz}{\field{C}}
\newcommand{\zz}{\field{Z}}
\newcommand{\nz}{\field{N}}
\newcommand{\R}{{\mathbb R}}
\newcommand{\Z}{{\mathbb Z}}
\newcommand{\I}{\mathrm{i}}
\newcommand{\Id}{\mathrm{Id}}

\renewcommand{\Re}{\operatorname{Re}}
\renewcommand{\Im}{\operatorname{Im}}

\newcommand{\cH}{\mathcal{H}}

\newcommand{\supp}{\mathrm{supp}}
\newcommand{\T}{\mathcal{T}}

\newcommand{\id}{1 \hspace{-1mm}  \mbox{\sf I}}

 \newtheorem{theorem}{Theorem}[section]

 \newtheorem{lem}[theorem]{Lemma}
 \newtheorem{cor}[theorem]{Corollary}
 \newtheorem{pro}[theorem]{Proposition}
 
 \newtheorem{rem}[theorem]{Remark}

\title[Scattering at low energies]{Scattering at low energies on manifolds 
with cylindrical ends and stable systoles}

\author[W. M\"uller]{Werner M\"uller}
\address{Mathematisches Institut, Universit\"at Bonn, Beringstrasse 1,
D-53115 Bonn, Germany} \email{mueller@math.uni-bonn.de}

\author[A. Strohmaier]{Alexander Strohmaier}
\address{Department of Mathematical Sciences,  Loughborough University,  Loughborough, Leicestershire, LE11 3TU,
UK} \email{a.strohmaier@lboro.ac.uk}

\begin{document}
\begin{abstract}
 Scattering theory for $p$-forms on manifolds with cylindrical ends has a 
direct interpretation in terms of cohomology. Using the Hodge isomorphism,
the scattering matrix at low energy may be regarded as operator on the
cohomology of the boundary. 
 Its value at zero describes the image of the absolute cohomology in the
 cohomology of the boundary. We show that the so-called scattering length, the
 Eisenbud-Wigner time delay at zero energy, has a cohomological
 interpretation as well. 
 Namely, it relates the norm of a cohomology class on the boundary to the
 norm of its image under the connecting homomorphism
 in the long exact sequence in cohomology. An interesting consequence of
 this is that one can estimate the scattering lengths in terms of geometric
 data like the volumes of certain homological systoles. 
\end{abstract}

\vspace{3cm}

\thanks{The second author was supported by the Leverhulm trust and the MPI Bonn}

\maketitle

\section{Introduction and main results}
\setcounter{equation}{0}

\begin{figure}[!h]
  \centerline{\includegraphics*[width=9cm]{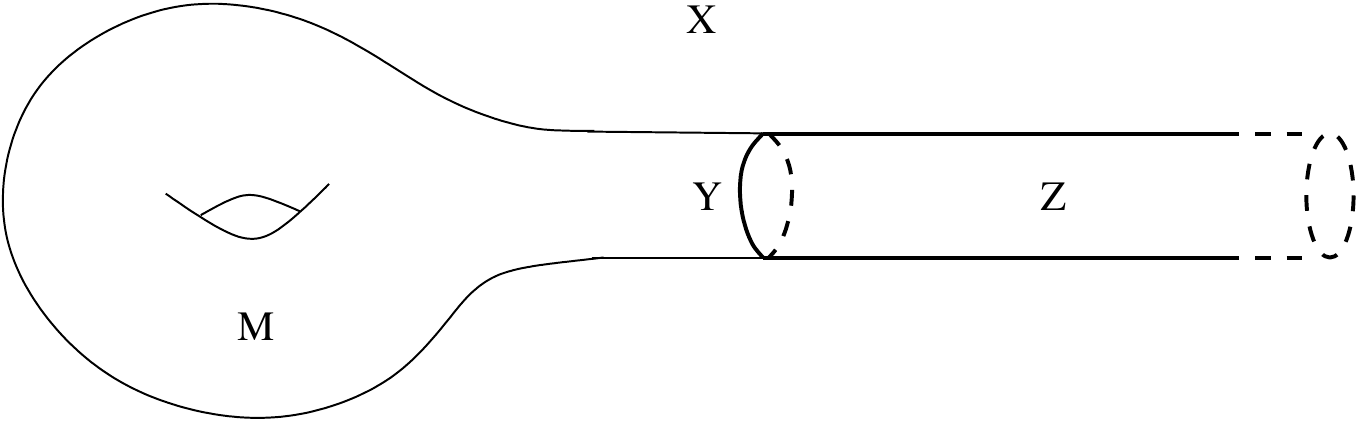}}
  \caption{Elongation $X$ of $M$.} \label{fiig01}
\end{figure}

Scattering theory for manifolds with cylindrical ends deals with the following
geometric situation. 
Let $M$ be an oriented, connected, compact Riemannian manifold with boundary
$Y=\partial M$ such that the metric is a product  near the boundary, i.e., 
there is a tubular neighborhood of $Y$ which is isometric to
$(-\epsilon,0] \times Y,$
equipped with the product metric $du^2 + h$, where $h$ is a Riemannian metric 
on $Y$.
The non-compact elongation $X$ of $M$ is then obtained from $M$ by
attaching the half-cylinder $Z=\rz^+\times Y$ over the boundary 
(see figure \ref{fiig01}):
\begin{equation}
 X=Z\cup_{Y} M.
\end{equation}

The Riemannian metric on $M$ is extended to one on $X$ in the 
obvious way so that
\begin{gather}
 g|_{\rz^+ \times Y}= du^2 + h.
\end{gather}
Such a manifold is  called a {\it manifold with cylindrical ends}.

Scattering theory on $X$ investigates how wave packets coming in from
infinity are scattered in $M$. The scattering of $p$-forms on $X$ 
is described by the scattering matrix
\begin{gather}
 C_p(\lambda) \in \mathrm{End}\left( \bigoplus_{\mu \leq \lambda}
   \mathrm{Eig}_\mu(\Delta'_p) \oplus \mathrm{Eig}_\mu(\Delta'_{p-1}) \right),
\end{gather}
where $\Delta'_p$ is the Laplace-Beltrami operator acting on $p$-forms on $Y$
and $\mathrm{Eig}_\mu(\Delta'_p)$ the eigenspace of $\Delta'_p$ with 
eigenvalue $\mu$.
In particular for small values of the spectral parameter $\lambda$ we have
\begin{gather}
 C_p(\lambda) \in \mathrm{End}\left( \mathcal{H}^p(Y) \oplus \mathcal{H}^{p-1}(Y) \right),
\end{gather}
where $\mathcal{H}^p(Y) = \ker \Delta_p'$ is the  space of harmonic $p$-forms
on $Y$.
For the purposes of this article one can think of the scattering matrix
for small values of the spectral parameter as being defined by the statement
of Theorem \ref{main1}. 
In this case it can be shown that $C(\lambda)$ leaves the direct summands
invariant and is of the form
\begin{gather}
  C_p(\lambda)=\left( \begin{matrix} S_p(\lambda) && 0\\ 0 && - S_{p-1}(\lambda) \end{matrix} \right),
\end{gather}
 where $S_p(\lambda) \in  \mathrm{End}( \mathcal{H}^p(Y))$ is the scattering matrix describing the scattering of coclosed
 forms. Again $S_p$ can be defined by the statement of Theorem \ref{main2} as
 the matrix relating the incoming and outgoing waves. The first
 observation is that the total scattering matrix for coclosed $p$-forms at
energy $0$, 
 \begin{gather*}
  S(0) = \oplus_p S_p(0) \in \mathrm{End}(\mathcal{H}^*(Y)),
 \end{gather*}
 is a self-adjoint involution which anti-commutes with the Hodge star operator.
 The $+1$ eigenspace of $S(0)$ coincides with the
 space of harmonic forms that represent cohomology classes in 
$\mathrm{Im}(r\colon H^*(X,\rz) \to H^*(Y,\rz))$. 

The Eisenbud-Wigner time-delay operator $\T_p(\lambda)$ describes the
time-delay a $p$-form-wave
undergoes when being scattered in $M$ (see Appendix \ref{app2}). It can be 
calculated from the
Eisenbud-Wigner formula (see Appendix \ref{app2}), and  for small $\lambda$ it 
is given by
\begin{gather*} 
 \T_p(\lambda)=- \I C_p(\lambda)^* \frac{d}{d \lambda} C_p(\lambda).
\end{gather*}
Of course
\begin{gather} 
  \T_p(\lambda)=\left( \begin{matrix} T_p(\lambda) && 0\\ 0 &&  T_{p-1}(\lambda) \end{matrix} \right),
 \end{gather}
 where $T_p(\lambda)$ is the time delay operator for coclosed forms defined
 by
 \begin{gather}
  T_p(\lambda)=-\I S_p(\lambda)^{-1} S'_{p}(\lambda).
 \end{gather}
Its value $T_p(0)$ at zero energy is of particular interest and we call it the
{\it scattering length}. The physical interpretation of the scattering length 
is as follows. If a 
coclosed wave packet has very low energy then, by the uncertainty relation, it 
will be far
spread out. In particular it will not be able to ``feel'' details of the
geometry of $M$. The effect of the manifold $M$ in the scattering process for 
a wave is then close to that of a cylindrical obstacle of length given by the 
scattering length.
It is therefore an interesting question to determine what geometric
properties of $M$ have an effect on the scattering length.
Since $T(0)$ commutes with the Hodge star operator it is enough
to know its restriction to the $-1$ eigenspace of $S(0)$. Denote by 
$\| \cdot \|_{st}$ the stable norm of a homology class, and 
by
$\|\cdot\|_\infty$ the comass norm on the cohomology groups (see \cite{Gro99}). 
Let $\nu_1>0$ be the smallest positive eigenvalue of $\Delta_p^\prime$. Let us define
\[
\mathrm{Vol}_*(M)=\mathrm{Vol}(M)+\frac{1}{\sqrt{\nu_1}}\mathrm{Vol}(Y).
\]
Furthermore, in section 
\ref{lzweiestimates} we introduce for each $n\in\nz$ and $0\le p\le n$
constants $C(n,p)>0$, which are related to the estimation of the comass norm on
$\Lambda^p\R^n$. They are equal to 1 for $p=0$ and $p=1$. 
One of our main results relates the scattering length to certains norms in 
homology
(Theorem \ref{main3}) and gives rise to
the following  estimation of the scattering length in terms of geometric
data.

\begin{theorem}\label{main00}
Let $0\le p\le n$. For
every $\phi$ in the $-1$-eigenspace of $ S_p(0)$ we have
\begin{gather*}
 \frac{1}{2} C(n,p+1)^{-1} \mathrm{Vol}_*(M)^{-1} \| [M] \cap 
\partial \phi\|_{st}^2 \leq  \langle \phi, T(0)^{-1} \phi \rangle 
\leq 
\frac{1}{2} C(n,p+1) \mathrm{Vol}(M) \|\partial \phi \|_\infty^2,
\end{gather*}
where $\partial: H^p(Y) \to H^{p+1}(M,Y)$ is the connecting homomorphism in the long exact sequence
in cohomology and $[M] \cap \partial \phi$ is the Poincare dual of the class $\partial \phi$.
\end{theorem}

As an example we treat the case when $Y$ has two connected
components $Y_1$ and $Y_2$ and $p=0$. In this case there is a canonical basis in $H^0(Y,\rz)$
with respect to which $T_0(0)$ has the form
\begin{gather}
 T_0(0) = \left ( \begin{matrix} t_1 & 0 \\ 0 & t_2 \end{matrix} \right), 
\end{gather}
so that
\begin{gather}
 t_1 = 2 \frac{\mathrm{Vol}(M)}{\mathrm{Vol}(Y)},\\
 C_2 \leq t_2 \leq C_1,
\end{gather}
and the constants $C_1$ and $C_2$ are given by
\begin{gather}
 C_1= 2 \mathrm{Vol}_*(M) \frac{\mathrm{Vol}(Y_1) \mathrm{Vol}(Y_2)}
 {\|\iota_*([Y_1])\|_{st}^2 (\mathrm{Vol}(Y_1)+ \mathrm{Vol}(Y_2))},\\
 C_2= 2 \mathrm{Vol}(M)^{-1} \frac{\mathrm{dist}(Y_1,Y_2)^2 \mathrm{Vol}(Y_1)\mathrm{Vol}(Y_2)}
 {\mathrm{Vol}(Y_1)+ \mathrm{Vol}(Y_2)}.
\end{gather}
The map $\iota$ is the inclusion of $Y$ into $M$. So we get an estimate for the
scattering length by purely geometric quantities.  
The physical interpretation of this is as follows. For a wave of low energy the reflection
coefficient $r_{11}$ and the transmission coefficient $r_{12}$ for a wave coming in at the boundary component
$Y_1$ are approximated by their values at zero, namely (see section \ref{scatterex1}) 
\begin{gather}
 r_{11}=\frac{\mathrm{Vol}(Y_1)-\mathrm{Vol}(Y_2)}{\mathrm{Vol}(Y_1)+\mathrm{Vol}(Y_2)}, \qquad r_{12}=\frac{2 \mathrm{Vol}(Y_1)}{\mathrm{Vol}(Y_1)+\mathrm{Vol}(Y_2)}.
\end{gather}
The time-delay is then determined by $t_1$ and $t_2$. For example in the case
where $\mathrm{Vol}(Y_1)=\mathrm{Vol}(Y_2)$ the reflection coefficient at zero 
energy
is zero and the time delay of the transmitted wave is equal to 
$\frac{1}{2}(t_1+t_2)$ 
(see section \ref{scatterex1}).
Another example, the full-torus, is treated in section \ref{scatterex2}.

Given $a>0$, let $M_a$ be the manifold that is 
obtained from $M$ by attaching the cylinder $[0,a] \times Y$ to $M$. 
We also investigate how the $L^2$-norm of a class in $H^p(M_a,\rz)$
and the $L^2$-norm of its image in $H^{p+1}(M_a,\rz)$ under the connecting
homomorphism are related for the manifold with boundary $M_a$. 
There is an operator $q_a$ that relates the
$L^2$-norm of classes in the complement of the kernel of the connecting
homomorphism to the $L^2$-norm of the image of this class under the connecting
homomorphism. So $q_a$ measures to what extent the connecting homomorphism 
is a partial isometry.
Theorem \ref{theoremquadratic} shows that the operator $q_a$ has an expansion
of the form
\begin{gather}
 q_a = a \mathbf{1} + \frac{1}{2} T(0) + O(e^{-c a}),
\end{gather}
as $a \to \infty$.
This shows that one can calculate the scattering length by approximating $X$
by the compact manifolds $M_a$ and consider the constant term in the above
expansion. The exponential decay of the remainder term may also be useful for 
numerical computations.

The paper is organized as follows. In sections \ref{one} we review
stationary scattering theory for $p$-forms on manifolds with cylindrical ends.
Section \ref{three} and section \ref{four} deal with cohomology of $M$ and
$X$, and their relation to scattering theory and the continuous spectrum of
the Laplacian of $X$. We also derive a cohomological formula for the 
scattering length. In sections \ref{lzweiestimates} and \ref{systoles} it is 
shown that the
$L^2$-scalar products on the cohomology groups of $X$ can be estimated in
terms of geometric quantities and that these estimates imply estimates on
the scattering length. Section \ref{examples} treats the case of functions
in the case of two boundary components and the case of a full-torus. In this
section we demonstrate how our main result can be used to obtain estimates of
the scattering length in terms of  geometric data.
In appendix
\ref{app2} we discuss the relation between the stationary and the dynamical 
approach to scattering theory 
and we establish the Eisenbud-Wigner formula for manifolds with cylindrical 
ends.

Whereas for the sake of notational simplicity we restricted ourselves in this paper to
manifolds with cylindrical ends most of our analysis carries over in a
straightforward manner to waveguides if Neumann boundary conditions are imposed.

\noindent{\bf Acknowledgements.}
Much of the work on this paper has been done during the visit of the second
author to the MPI in Bonn and he would like to thank the MPI for the kind
support. Both authors would also like to
thank the MSRI in Berkeley for hospitality during the
program ``Analysis on Singular Spaces''. We are grateful to Werner Ballmann, Alexej Bolsinov,
Peter Perry, and Sasha Pushnitski for useful discussions and comments.
The authors would like to thank the referee for his very useful comments and suggestions.

\section{Stationary scattering theory for manifolds  with cylindrical 
ends} \label{one}

In this section we review stationary scattering theory for differential forms on manifolds with cylindrical ends
and establish some elementary relations which we will use in subsequent sections.
For functions this can be found e.g in \cite{Chr95,Chr02,Mel93} or in \cite{IKL10} which also covers the asymptotically cylindrical case.
Scattering theory for Dirac type operators that are of product type on the cylinder is also well covered in the literature (see e.g.
\cite{Mu94})

As before let $M$ be an oriented, connected, compact Riemannian manifold with
boundary $Y=\partial M$ such that the metric is a product  near the boundary.
Let $X$ be the elongation of $M$. For any Riemannian manifold $W$ we
 denote by $\Lambda^p(W)$ (resp. 
$\Lambda_c^p(W)$, $L^2\Lambda^p(W)$) the space of
smooth $p$-forms  (resp. smooth $p$-forms with compact supports, $L^2$-forms) 
on $W$. Let $\Delta_p$ be the Laplace operator on $\Lambda^p(W)$. 
Throughout this paper a harmonic $p$-form will mean a $p$-form 
$\phi\in\Lambda^p(W)$ with $\Delta_p\phi=0$.

Since $X$ is complete, the Laplace-Beltrami operator $\Delta_p$ on $p$-forms is
essentially self-adjoint when regarded as operator in $L^2\Lambda^p(X)$ with
domain $\Lambda_c^p(X)$ (\cite{Che73}). We continue to denote its self-adjoint
extension by $\Delta_p$. In this section we recall some facts concerning the
generalized eigenforms of $\Delta_p$ and derive some properties of the 
scattering matrix for low energy.

Let $\Delta'_p$ denote the Laplacian on $p$-Forms of $Y$. Let 
$0\le\nu_1<\nu_2<\cdots$ be the distinct eigenvalues of 
 of $\Delta'_p\oplus\Delta'_{p-1}$. Let $\Sigma\to \cz$
the minimal Riemann surface on which $\sqrt{\lambda^2-\nu_j}$ is a 
single-valued
function for all $j\in\nz_0$. As proved by Melrose\cite{Mel93}  the resolvent 
$(\Delta_p-\lambda^2)^{-1}$, regarded as operator $\Lambda_c^p(X)\to 
L^2_{\mathrm{loc}}\Lambda^p(X)$, admits a meromorphic extension from the 
half-plane
$\Im(\lambda)>0$ to $\Sigma$. Let $\mu_1^2>0$ be the first non-zero
eigenvalue
of $\Delta'_p\oplus\Delta'_{p-1}$. Then it follows in particular that
$(\Delta_p-\lambda^2)^{-1}$ extends to a meromorphic function on the disc
$\{\lambda \colon |\lambda|< \mu_1\}$. As a consequence one can define 
analytic families of generalized eigenforms.

For $a \geq 0$ we denote by  $Y_a$ the hypersurface $(a,Y) \subset 
\rz^+_0 \times Y\subset X$. 
Note that the restriction of the bundle $\Lambda^p T^* X$ to $Y_a$ 
is canonically isomorphic to the direct sum $\Lambda^p T^* Y \oplus
\Lambda^{p-1} T^* Y$ since each vector $f \in \Lambda^p T^* X$
at some point $(u,x)$ can be uniquely decomposed as
$f=f_1 + du \wedge f_2$ with $f_1 \in \Lambda^p T^* X$ and $f_2 \in
\Lambda^{p-1} T^* X$. 
Accordingly, the restriction of any $p$-form $\omega \in\Lambda^p(X)$
to the cylinder $\rz^+ \times Y$ is of the form
\begin{gather}
 \omega= \omega_1 + du \wedge \omega_2,
\end{gather}
where $\omega_1$ and $\omega_2$ are sections of the pulled back bundles 
$\pi^*\Lambda^p T^*
Y$ and $\pi^*\Lambda^{p-1} T^* Y$, respectively, and  
$\pi:  \rz^+ \times Y\to Y$ is the canonical projection.
We think of $\omega_1(u)$ and $\omega_2(u)$ as forms on $Y$
that depend smoothly on the additional parameter $u$. The map
\begin{gather}
 j_p: \pi^* \Lambda^p T^* Y \oplus \pi^* \Lambda^{p-1} T^* Y \to \Lambda^p T^* Z,
 \quad (\omega_1,\omega_2) \mapsto \omega_1 + du \wedge \omega_2
\end{gather}
is an isomorphism of vector bundles.
The exterior differential of such a form $\omega$ is then given by
\begin{gather}
 d \omega = d' \omega_1 + du \wedge \partial_u \omega_1 - du \wedge d'\omega_2,
\end{gather}
where $d'$ denotes the exterior differential on $Y$.
In matrix notation this means
\begin{gather}\label{differ}
 j^{-1} \circ d \circ j =\left( \begin{matrix} d' && 0 \\ \partial_u && -d' \end{matrix} \right),
\end{gather}
where we denote $j= \oplus j_p$.
Since the metric has product structure the decomposition is orthogonal and 
the formal adjoint $\delta$ of $d$ is therefore
\begin{gather}\label{codiffer}
 j^{-1} \circ \delta \circ j =\left( \begin{matrix} \delta' && -\partial_u \\ 0 && -\delta' \end{matrix} \right).
\end{gather}
Here again we use the notation $\delta'$ for the codifferential on
$Y$.
We therefore have
\begin{gather}
 j^{-1} \circ (d+\delta) \circ j=\left( \begin{matrix}  1 && 0\\ 0 && -1 
\end{matrix} \right)
 (d'+\delta') + \left( \begin{matrix}  0 && -1\\ 1 && 0 \end{matrix} \right) 
\partial_u.
\end{gather}
and
\[
j^{-1}_p \circ \Delta_p \circ j_p=
\left( \begin{matrix} -\partial_u^2+\Delta'_p  && 0\\ 0 && -\partial_u^2+
\Delta'_{p-1}\end{matrix} \right)
\]
for all $p$.
This form of an operator allows for a separation of variables on the cylinder
$\rz^+\times Y $. Suppose that $(\psi_i)$ is an orthonormal sequence of 
eigenforms of 
\begin{gather}
 \left( \begin{matrix}  \Delta_p' && 0\\ 0 && \Delta_{p-1}' \end{matrix}
 \right)
\end{gather}
with eigenvalues $\mu_{\psi_i}^2$. Then for $|\lambda|<\mu_1$,
any solution $F$ of the equation 
$(\Delta_p -\lambda^2) F=0$  has an expansion of the form
\begin{equation}\label{expansion1}
  F(u)=
\sum_{i=0}^\infty\left(a_i e^{-\I \sqrt{\lambda^2-\mu_{\psi_i}^2} u}  +  
b_i e^{\I \sqrt{\lambda^2-\mu_{\psi_i}^2} u}  \right)j_p(\psi_i)+\begin{cases}
\sum_{\mu_{\psi_i}=0}c_i uj_p(\psi_i),& \lambda=0;\\
0,&\lambda\not=0.
\end{cases}
\end{equation}
The series converges in the $C^\infty$-topology. The square roots are chosen throughout the article as
$\sqrt{r e^{\I \varphi}}=\sqrt{r} e^{\frac{\I}{2} \varphi}$ if $r>0$ and $0 \leq \varphi < 2 \pi$.
From the analytic continuation of the resolvent one gets the following result.

\begin{theorem} \label{main1}
 For each $\psi$ in $\mathrm{ker} (\Delta_p' \oplus \Delta_{p-1}')$
 there exists a $p$-form $\widetilde F(\psi,\lambda)$ which is meromorphic in 
$\lambda\in\{z \colon |z| < \mu_1\}$ such that the following conditions hold
 \begin{enumerate}
  \item[(i)] $\widetilde F(\psi,\lambda)$ is holomorphic in $\lambda$ for 
$\Im(\lambda)>0$.
  \item[(ii)] $\Delta \widetilde F(\psi,\lambda) = \lambda^2
\widetilde  F(\psi,\lambda)$.
  \item[(iii)] There exists $\widetilde R_p(\psi,\lambda)\in L^2\Lambda^p(Z)$ 
such that on $Z$ we have
\[
\widetilde F(\psi,\lambda)=e^{-\I \lambda u} j_p(\psi) + e^{\I \lambda u} 
j_p(C_p(\lambda)\psi) + \widetilde R_p(\psi,\lambda).
\]
\item[(iv)] $C_p(\lambda): \ker(\Delta_p' \oplus \Delta_{p-1}') \to
    \ker(\Delta_p' \oplus \Delta_{p-1}')$ is a linear operator, and 
$C_p(\lambda)$ and $\widetilde R_p(\psi,\lambda)$ are meromorphic functions of 
$\lambda$.
\end{enumerate}
 Moreover, $C_p(\lambda)$, $\widetilde R_p(\psi,\lambda)$
 and $\widetilde F(\psi,\lambda)$ are uniquely determined by these properties.
\end{theorem}
\begin{proof} This follows from\cite{Gui89},  \cite{Mel93}. For the 
convenience of the 
reader we include the details. Let $\Sigma\to\cz$ be the minimal 
Riemann surface to which all the functions $(\lambda^2-\mu^2_{\psi_i})^{1/2}$
extend to be holomorphic. 
By \cite[Th\'eor\`em 0.2]{Gui89}, \cite[Theorem7.1]{Mel95} the resolvent 
$(\Delta_p - \lambda^2)^{-1}$, regarded as operator 
$\Lambda^p_c(X)\to L^2_{\mathrm{loc}}\Lambda^p(X)$, extends to a meromorphic
function of $\lambda\in\Sigma$. Especially
$(\Delta_p - \lambda^2)^{-1}$ extends to a meromorphic function of $\lambda\in
\{z\in\cz\colon |z|<\mu_1\}$ as an operator
 $\Lambda^p_c(X)\to L^2_{\mathrm{loc}}\Lambda^p(X)$. 
 Let $\chi$ be a smooth function with support in
 $Z$ which is equal to $1$ outside a compact set. Put 
 \begin{gather}
  \widetilde F(\psi,\lambda):= \chi e^{-\I \lambda u} j_p(\psi) -  
(\Delta_p - \lambda^2)^{-1}\left\{ (\Delta_p - \lambda^2)
 (\chi e^{-\I \lambda u} j_p(\psi))\right\}.
 \end{gather}
 Then $\widetilde F(\psi,\lambda)$ is a smooth $p$-form on $X$ which depends 
meromorphically on $\lambda\in\Sigma$. Moreover, it satisfies
 \begin{gather}
  (\Delta_p - \lambda^2) \widetilde F(\psi,\lambda)=0.
 \end{gather}
 Since $\widetilde F(\psi,\lambda)-\chi e^{-\I \lambda u} j_p(\psi)$ is square 
integrable 
for $\Im(\lambda) > 0$, the expansion (\ref{expansion1})  has the form
\[
\widetilde F(\psi,\lambda)=e^{-\I \lambda u} j_p(\psi) + e^{\I \lambda u} 
j_p(C_p(\lambda)\psi) + \widetilde R_p(\psi,\lambda)
\]
where $C_p(\lambda)\psi\in\ker(\Delta'_P\oplus\Delta'_{p-1})$ and 
\begin{gather} \label{expansion}
   \widetilde R_p(\psi,\lambda,u)=
   \sum_{\mu_{\psi_i}\not=0} 
  \left(a_i(\lambda) e^{-\I \sqrt{\lambda^2-\mu_{\psi_i}^2} u} j_p(\psi_i) +  
b_i(\lambda)e^{+\I \sqrt{\lambda^2-\mu_{\psi_i}^2} u} j_p(\psi_i) \right).
\end{gather}
Moreover, $\widetilde R_p(\psi,\lambda)$ is square integrable for 
$\Im(\lambda)>0$. 
This implies that $a_i(\lambda)=0$ for $\Im(\lambda)>0$. Since $a_i(\lambda)$
is meromorphic, it follows that $a_i(\lambda)=0$ for $|\lambda|<\mu_1$. 
Thus we conclude that $R_p(\psi,\lambda)$ is a meromorphic function in the
disc $|\lambda|<\mu_1$ with values in $L^2\Lambda^p(X)$.

The uniqueness is an immediate consequence of the self-adjointness 
of $\Delta_p$.
 Namely, if $\widetilde F_1(\psi,\lambda)$ and $\widetilde F_2(\psi,\lambda)$ 
both have the above properties
then their difference $G(\psi,\lambda)$ is square integrable for 
$\Im(\lambda)>0$
  and is contained in the kernel of $(\Delta_p-\lambda^2)$. 
Using that $\Delta_p$ is self-adjoint, we get $G(\psi,\lambda)=0$ for 
$\Im(\lambda)>0$.
Since it is meromorphic in $\lambda$ we conclude that $G(\psi,\lambda)=0$. 
\end{proof}
The requirement that $\widetilde R_p(\psi,\lambda)$ is in $L^2\Lambda^p(Z)$
in fact implies a much faster decay at infinity.
\begin{lem}\label{rapiddec}
If $\widetilde R(\psi,\lambda)$ is regular at $\lambda$, then 
all  derivatives of $\widetilde R(\psi,\lambda)$ in $u$ and $x$
are exponentially decaying as $u \to \infty$. More precisely
\[
  |\partial_u^k (\Delta')^{l} \widetilde R(\psi,\lambda) | \leq C_{k,l,\lambda}' 
e^{-C_\lambda u},
\]
 for all $k,l \geq 0$ and some positive constants $C_{k,l,\lambda}'$ and 
$C_\lambda$.
\end{lem}
\begin{proof}
 The proof is already implicitly contained in the proof of the previous 
theorem.
 Namely, $\widetilde R(\psi,\lambda)$ is smooth by elliptic regularity and 
thus the expansion (\ref{expansion}) converges in $\Lambda^p(Z)$.   
Since  $a_i=0$, 
 we get the decay for $|\lambda| < \mu_1$ with 
$C_\lambda=\Re{\left(\sqrt{\mu_1^2-\lambda^2}\right)}$.
\end{proof}

Let $\psi\in\ker\Delta_p^\prime$. Put
\begin{equation}\label{eigenfunct1}
F(\psi,\lambda):=\widetilde F\left((\psi,0),\lambda\right),\quad 
F(du\wedge \psi,\lambda):= \widetilde F\left((0,\psi),\lambda\right).
\end{equation}
Then the expansion (iii) of Theorem \ref{main1} takes the form
\begin{equation}\label{expansion11}
\begin{split}
F(\psi,\lambda)&=e^{-i\lambda u}\psi+e^{i\lambda u}C_p^{11}(\lambda)\psi
+e^{i\lambda u}du\wedge C_p^{21}(\lambda)\psi+R_p(\psi,\lambda);\\
F(du\wedge \psi,\lambda)&=e^{-i\lambda u}du\wedge \psi
+e^{i\lambda u}du\wedge C_p^{22}(\lambda)\psi
+e^{i\lambda u} C_p^{12}(\lambda)\psi+R_p(du\wedge \psi,\lambda),
\end{split}
\end{equation}
where
\begin{equation}\label{scattmatrix}
C_p(\lambda)=\begin{pmatrix} C_p^{11}(\lambda)&C_p^{12}(\lambda)\\
{}&{}\\
C_p^{21}(\lambda)&C_p^{22}(\lambda)
\end{pmatrix}
\end{equation}
as endomorphism of $ \mathrm{ker} \Delta_p' \oplus \mathrm{ker}\Delta_{p-1}'$.
A priory there is  no reason why $C_p(\lambda)$ should leave 
the summands invariant. Nevertheless this is guaranteed by a continuous 
version of the Hodge decomposition as the proof of the following proposition 
shows.
\begin{pro}\label{vanishing}
 $C_p(\lambda)$ leaves the spaces $\mathrm{ker }\Delta_p'$ and $\mathrm{ker }
\Delta_{p-1}'$ invariant, i.e., $C_p^{12}(\lambda)=0$ and $C_p^{21}(\lambda)=0$.
\end{pro}
\begin{proof}
 
 Let $\psi \in \Delta_p'$.
Using the expansion of  \ref{expansion11}, it follows that on $\rz^+ \times Y$
we have
\begin{equation}\label{gl234}
\begin{split}
\delta d F(\psi,\lambda) &= 
\lambda^{2}e^{-\I \lambda u} \psi + \lambda^{2}
e^{\I \lambda u}  C_p^{11}(\lambda)\psi +   
\delta d R_p(\psi,\lambda),\\
d \delta F(du\wedge \psi,\lambda) &= 
\lambda^{2}e^{-\I \lambda u} du\wedge \psi + 
\lambda^{2}e^{\I\lambda u} du\wedge  C_p^{22}(\lambda) \psi + d \delta 
R_p(du\wedge \psi,\lambda).
\end{split} 
\end{equation}

Comparing the leading terms,  it follows  that 
\[
\lambda^{-2} \delta d F(\psi,\lambda)=F(\psi,\lambda),\quad 
\lambda^{-2} d \delta F(du\wedge \psi,\lambda)=F(du\wedge \psi,\lambda).
\]
 Since all derivatives of $R_p(\psi,\lambda)$ and $R_p(du\wedge \psi,\lambda)$
are exponentially decaying,
 the uniqueness statement in theorem \ref{main1} implies immediately
 $C_p^{12}(\lambda)=0$ and $C_p^{21}(\lambda)=0$.
\end{proof}

Note that the splitting $\Lambda^p T^*Z=\pi^*(\Lambda^p T^*Y) \oplus  
\pi^*(\Lambda^{p-1}T^*Y)$ indeed
corresponds to the Hodge decomposition. Let $P_1$ and $P_2$ the projection on
the first and second summand, respectively.  Then we have
\begin{pro}
 We have $P_2 \psi=0$, i.e. $\psi \in \mathrm{ker} \Delta_p'$ if and only if
 $\delta \widetilde F(\psi,\lambda)=0$. Similarly $P_1 \psi=0$, i.e. 
$\psi \in \mathrm{ker} \Delta_{p-1}'$ if and only if $d \widetilde 
F(\psi,\lambda)=0$. 
\end{pro}
\begin{proof}
 Let $\psi=(\psi_1,\psi_2)\in\ker\Delta_p^\prime\oplus\ker\Delta_{p-1}^\prime$.
Suppose that $\delta \widetilde F(\psi,\lambda)=0$. Then, in particular, 
we have $\delta F(du\wedge \psi_2,\lambda)=0$. Applying $\delta$ to the
second equation of (\ref{expansion11}) and using Proposition \ref{vanishing},
 it follows  that $\psi_2=0$.
 For the other direction, observe that by (\ref{expansion11}),
$\delta F(\psi_1,\lambda)=
\delta R(\psi_1,\lambda)$ on $Z$. Hence, $\delta F(\psi_1,\lambda)$ is 
exponentially decaying. 
 In particular it is square integrable for $\Im (\lambda) >0$.
 Since $\Delta_p$ is self-adjoint and
 $(\Delta_p-\lambda^2) \delta F(\psi_1,\lambda)=0$ we get 
$\delta F(\psi_1,\lambda)=0$. Thus, if $\psi_2=0$, it follows that 
$\delta\widetilde F(\psi,\lambda)=0$. 
 The proof of the other case is analogous.
\end{proof}

\begin{pro}
 The following relation holds for $0 \leq p < n$.
  \begin{gather}
 C_p(\lambda) |_{\mathrm{ker} \Delta_p'} = -C_{p+1}(\lambda) |_{\mathrm{ker} \Delta_p'},
 \end{gather}
\end{pro}
\begin{proof} 
Let $\psi \in \mathrm{ker} \Delta_p'$. Using (\ref{expansion11}), we get
\begin{gather}
\I \lambda^{-1} d F(\psi,\lambda) |_Z=  e^{-\I \lambda u} du \wedge \psi 
    - e^{\I \lambda u} du \wedge C_p^{11}(\lambda) \psi +   
\I \lambda^{-1}d R_p(\psi,\lambda)
\end{gather}
 Comparing the leading terms, it follows from Theorem \ref{main1} that
\[
\I \lambda^{-1} d F(\psi,\lambda)=F(du\wedge\psi,\lambda).
\]
Using (\ref{expansion11}), we get $C_p^{11}(\lambda)=-C_{p+1}^{22}(\lambda)$,
which is equivalent to the statement of the Proposition.
\end{proof}
 
 Let us use the notation $S_p(\lambda)$ for the restriction of $C_p(\lambda)$
 to $\mathrm{ker} \Delta_p'$. Then the above proposition shows that
 \begin{gather} \label{invsummands}
  C_p(\lambda)=\left( \begin{matrix} S_p(\lambda) && 0\\ 0 && - S_{p-1}(\lambda) \end{matrix} \right).
 \end{gather}
 
In the following we will  suppress the index $p$ and write $S(\lambda)$, 
meaning that $S(\lambda)$ is acting on the space of harmonic forms, 
leaving the space of $p$-forms invariant. It is the scattering matrix at low 
energy for the scattering problem for coclosed forms.
 Summarizing we have established the following theorem.
 
\begin{theorem} \label{main2}
 For each harmonic $p$-form $\psi \in \mathrm{ker} (\Delta_p')$
 there exists a $p$-form $F(\psi,\lambda)$ on $X$, which is meromorphic in 
$\lambda$ in the disc $|\lambda| < \mu_1$ such that
 \begin{enumerate}
  \item[(i)] $\delta F(\psi,\lambda)=0$
  \item[(ii)] $F(\psi,\lambda)$ is holomorphic in $\lambda$ for $\Im(\lambda)>0$.
  \item[(iii)] $(\Delta_p - \lambda^2) F(\psi,\lambda)=0$.
  \item[(iv)] $F(\psi,\lambda)=
         e^{-\I \lambda u} \psi + e^{\I \lambda u} S(\lambda)
         \psi + R(\psi,\lambda)$ on $\rz^+ \times Y$,
  \item[(v)] $S(\lambda) \in \mathrm{End}(\ker(\Delta_p'))$ and
    $R(\psi,\lambda) \in L^2\Lambda^p( Z)$
 are meromorphic functions of  $\lambda$.
 \end{enumerate}
 Moreover, $S(\lambda)$,
 $R(\psi,\lambda)$
 and $F(\psi,\lambda)$ are uniquely determined by these properties.
\end{theorem}
 
The scattering matrix has the following properties
\begin{theorem}\label{main22}
 The function $S(\lambda)$ satisfies the following equations
 \begin{enumerate}
  \item[(i)] $S(\lambda)^* S(\overline \lambda)=1$.
  \item[(ii)] $S(\lambda) S(-\lambda)=1$.
  \item[(iii)] $S(\lambda) * = - * S(\lambda)$,
 \end{enumerate}
 where $*$ is the Hodge star operator on $Y$.
\end{theorem}
\begin{proof}
 For $a>0$ let $X_a$ be the manifold $([0,a)\times Y) \cup_Y X$
 with boundary $Y_a=\{a\} \times Y$. Let $\omega_X$ be the volume form of $X$. 
By Theorem \ref{main2}, (iii), we have
 \begin{gather} \nonumber
  \int_{X_a} \langle F(\psi,\lambda),\Delta_p
F(\psi,\overline\lambda)
  \rangle \omega_X 
  - \int_{X_a} \langle \Delta_p F(\psi,\lambda),
 F(\psi,\overline\lambda) \rangle \omega_X = 0.
 \end{gather}
Using Green's formula, we obtain
 \begin{gather}
  \int_{Y_a} \langle F(\psi, \lambda),-\partial_u F(\psi,\overline\lambda)
  \rangle \omega_Y
  + \int_{Y_a} \langle \partial_u F(\psi, \lambda),F(\psi,\overline \lambda) \rangle \omega_Y= 0.
 \end{gather}
 In the limit $a \to \infty$ this expression can be evaluated using Theorem
 \ref{main2}. We obtain
 \begin{gather} \nonumber
  \lim_{a\to \infty}\int_{Y_a} ( \langle F(\psi,\lambda),-\partial_u 
F(\psi,\overline \lambda)
  \rangle
  + \langle \partial_u F(\psi,\lambda),F(\psi,\overline \lambda) \rangle) \omega_Y =\\
  =- 2 \I \lambda (\Vert \psi \Vert^2-\langle S(\lambda)
  \psi,S(\overline\lambda) \psi\rangle)= 0,
 \end{gather}
 which proves the first statement.
 The second statement follows from the functional equation
\[
 F(C(\lambda)\psi,-\lambda)=F(\psi,\lambda),
\]
 which is a simple consequence
 of the uniqueness statement in theorem \ref{main1}. To show that
 $C(\lambda)$ anticommutes with the Hodge star operator on $Y$ we note that
 that Hodge star operator $*_X$ on $X$ commutes with the Laplace operator 
$\Delta$, i.e.
 $*_X \Delta_p = \Delta_{n-p} *_X$. Applying this to $F(\psi,\lambda)$ and 
using the uniqueness
 statement we obtain immediately
 $$
  *_X C_p =  C_{n-p} *_X.
 $$
  For $\psi \in \mathrm{ker} \Delta_p'$ we get $*_X \psi = (-1)^p du 
\wedge * \psi$
 and consequently $* S(\lambda)= -S(\lambda) *$, where we used that $du \wedge$
 anticommutes with $C(\lambda)$.
\end{proof}

As an application we obtain the following well known result about the signature
$\mathrm{sign}(Y)$ of a closed manifold $Y$.

\begin{cor}
Let $Y$ be a closed oriented manifold. Assume that $Y$ is the boundary of a 
compact manifold. Then $\mathrm{sign}(Y)=0$.
\end{cor}
\begin{proof} 
We may assume that $\dim Y=4k$. Otherwise the signature is zero. Pick a 
Riemannian metric on $Y$. Let $\cH^{2k}_\pm(Y)$ be the 
$\pm1$ eigenspaces of $*$ acting in $\cH^{2k}(Y)$. Then the signature 
$\mathrm{sign}(Y)$ of $Y$ is given by
\[
\mathrm{sign}(Y)=\dim\cH^{2k}_+(Y)-\dim\cH^{2k}_-(Y).
\]
If $Y$ is the boundary of a compact Riemannian manifold $M$, it follows from
Theorem \ref{main22}, that $S(\lambda)$ is regular at $\lambda=0$, $S(0)^2=\id$
and $S(0)$ intertwines $\cH^{2k}_+(Y)$ and $\cH^{2k}_-(Y)$. Hence we get
$\mathrm{sign}(Y)=0$.
\end{proof}
\begin{rem}
The same proof works equally well for Dirac type operators. It
implies the cobordism invariance of the index of Dirac operators.
\end{rem}

\begin{pro} \label{reg}
 $S(\lambda)$, $F(\psi,\lambda)$, and $R(\psi,\lambda)$ are regular for real 
$\lambda$.
 If $\psi = -S(0)\psi$ then $F(\psi,0)=0$.
\end{pro}
\begin{proof}
 For real $\lambda$ it follows from Theorem \ref{main22} that
 $S(\lambda)S^*(\lambda)=\id$ and therefore 
$\|S(\lambda)\|=1$.
 In particular, $S$ is bounded on the real line and can not have a pole there.
 It remains to show that $F$ is regular for real $\lambda$. Suppose that
 $\phi$ is a square integrable eigensection of $\Delta$ with real eigenvalue
 $\lambda'$. Then the expansion (\ref{expansion1}) of $\phi$ on
$\rz^+\times Y$ takes the form
\[
\phi(u,y)=\sum_{\mu_{\psi_i}^2>\lambda^\prime}a_i
e^{-\sqrt{\mu_{\psi_i}^2-\lambda^\prime}\,u}j_p(\psi_i).
\]
This implies that $\phi$ is exponentially decaying.
 Since for real $\lambda$ 
 \begin{gather}
   0=\langle  (\Delta-\lambda^2) F(\psi,\lambda),\phi \rangle
=-(\lambda^2-\lambda'^2) \langle  F(\psi,\lambda),\phi \rangle, 
 \end{gather} 
 and $\langle  F(\psi,\lambda),\phi \rangle$ is a meromorphic function we get
 $ \langle  F(\psi,\lambda),\phi \rangle=0$.
 Suppose now that $F(\psi,\lambda)$ has a pole of order $k$ at $\lambda'$.
 Then $G:=\lim_{\lambda \to \lambda'} (\lambda-\lambda^\prime)^k
 F(\psi,\lambda)$ is an eigenform with eigenvalue $\lambda'$
 and it also is square integrable since $S$ is regular at $\lambda'$.
 By the above $\langle G, G \rangle=0$. It follows
 that $G=0$ and therefore $F(\psi,\lambda)$ is regular at $\lambda'$.
 If  $\psi = -S(0)\psi$ then $F(\psi,0)$ is square integrable and harmonic and 
by the same argument $F(\psi,0)=0$.
\end{proof}

By the above $F(\psi,0)$ and $F^\prime(\psi,0):=\frac{\partial}{\partial \lambda}
F(\psi,\lambda)|_{\lambda=0}$ are well defined. From the proof of Proposition
\ref{reg} we obtain the following corollary.

\begin{cor}\label{orthogonal}
 $F(\psi,0)$ and $F^\prime(\psi,0)$ are orthogonal to the space
$\cH^p_{(2)}(X)$ of square integrable harmonic forms.
\end{cor}

The scattering matrix $S(\lambda)$ is also regular at $0$ and it follows 
from Theorem \ref{main22} that $S(0)$
is a self-adjoint involution. Hence, $\mathrm{ker}(\Delta_p')$ decomposes
into $+1$ and $-1$ eigenspaces.
For $\psi \in \mathrm{ker}(\Delta_p')$ with $S(0)\psi=\psi$ we get
that $F(\psi,0)$ is a smooth coclosed harmonic $p$-form whose restriction
to $Z$ equals 
$$
F(\psi,0)= 2\psi+ R(\psi,0), 
$$
where $R(\psi,0)$ and its derivatives are exponentially decaying.
That is, $\psi$ is a limiting value of $\frac{1}{2}F(\psi,0)$
in the sense of Atiyah, Patodi, and Singer (\cite{APS75}).
It turns out that the converse is also true.
\begin{pro}\label{limiting}
 The $+1$ eigenspace for $S(0)$ is the set of limiting values of
  coclosed harmonic forms on $X$, i.e, it equals 
\[
\{\psi \in \mathrm{ker} \Delta_p' \mid \exists G \in 
\Lambda^p(X)\colon
G|_Z - \psi \in L^2\Lambda^p(Z), \; \Delta_p G =0,\;\delta G=0\}.
\]
Furthermore for each $\psi\in\ker \Delta_p^\prime$, $F(\psi,0)$ satisfies
$d F(\psi,0)=0$ and $\delta F(\psi,0)=0$. 
\end{pro}
\begin{proof}
 Suppose that $F \in \Lambda^p(X)$ and 
 $G \in \Lambda^{n-1-p} (X)$ are both coclosed harmonic forms on 
$X$ with limiting values $\psi$ and $\phi$, respectively.
 Since $\psi-G$ and $\phi-F$ are both exponentially decaying
 and $\psi$ and $\phi$ are closed and coclosed on $Y$ the forms 
 $dG$ and $dF$ are exponentially decaying.
 Using Green's formula, we get
\[
0=\langle\Delta G,G\rangle=
\langle\delta dG,G\rangle =\langle dG,dG \rangle.
\]
Thus $dG=0$. Similarly we get $dF=0$.
Using Stokes formula, it follows that
 $$
  0=\int_{X_a} d G \wedge F= \pm \int_{Y} \psi \wedge \phi + O(e^{-cu}).
 $$
 Thus, $\langle \psi, *\phi \rangle=0$.
 Now suppose that $\phi$ is a limiting value which is in the $-1$
 eigenspace to $S(0)$. Since $*$ anticommutes with $S(0)$, it follows that
 $*\phi$ is in the $+1$ eigenspace. It is therefore a limiting value.
 Since $*\phi$ and $\phi$ are both limiting values it follows from the above
 that $\|\phi\|^2=\pm \langle \phi, **\phi\rangle=0$ and therefore, $\phi=0$.
 Since the set of limiting vectors contains the $+1$ eigenspace it has to
 coincide with the $+1$ eigenspace.
\end{proof}

Finally we derive some formulas concerning $F^\prime(\psi,0)$ which we are going
to use in the next section. Note that the restriction of 
 $F'(\psi,\lambda)$ to the cylinder $Z$ has the form
\begin{gather}
 F'(\psi,\lambda)|_Z = -\I u (e^{-\I  \lambda u}\psi-e^{+\I  \lambda u}
S(\lambda) \psi) +e^{\I  \lambda u} S'(\lambda)\psi + R'(\psi,\lambda),
\end{gather}
and for $\lambda=0$:
\begin{equation}\label{expansion2}
 F'(\psi,0)|_Z = - \I u \left(1-S(0)\right) \psi+
 S'(0) \psi + R'(\psi,0).
\end{equation}

Differentiating the equation
\begin{gather}
 (\Delta - \lambda^2) F(\psi,\lambda)=0
\end{gather}
it follows that
\begin{gather}
 \Delta F'(\psi,0)=0.
\end{gather}

Hence, $dF'(\psi,0)$ is in the kernel of $\Delta$ and its restriction to the cylinder has the form
\begin{gather}
dF'(\psi,0)|_Z= -\I du \wedge (1-S(0)) \psi + dR'(\psi,0).
\end{gather} 
By Theorem \ref{main22}, (iii), we get
\begin{equation}\label{limitvalue}
*_X dF'(\psi,0)|_Z= -\I  (1+S(0)) * \psi + *_X dR'(\psi,0).
\end{equation} 
Thus, $*_X dF'(\psi,0)$ is an extended harmonic form with limiting value
$-\I  (1+S(0)) * \psi$.
Since $*_X dF'(\psi,0)$ is coexact and bounded it is orthogonal
to the space $\mathcal{H}_{(2)}^*(X)$ of square integrable harmonic forms.
The proof of Prop. \ref{reg} shows that $\frac{1}{2}F(-\I  (1+S(0)) * \psi,0)$
is harmonic and orthogonal to $\mathcal{H}_{(2)}^*(X)$. Their difference is 
therefore square integrable, harmonic and orthogonal to 
$\mathcal{H}_{(2)}^*(X)$.
Thus, it vanishes and we have the following nice formula
\begin{gather}\label{niceform}
*_X dF'(\psi,0) = -\frac{\I}{2} F\left((1+S(0)) * \psi,0\right).
\end{gather} 

In particular if $\psi$ is in the $-1$ eigenspace of $S(0)$ we have 
\begin{gather}
 *_X dF'(\psi,0) = - \I F(* \psi,0),
\end{gather}
or equivalently
\begin{equation}\label{fprime}
 dF'(\psi,0) = - \I F(du \wedge \psi,0).
\end{equation}

\begin{rem}
Note that whereas $F(\psi,\lambda)$ is regular at $\lambda=0$ the
meromorphic continuation of the resolvent $(\Delta_p-\lambda^2)^{-1}$  is in general not regular at zero
but has a pole of the form
$-K_1 \lambda^{-2} + (K_2 +K_3) \lambda^{-1}$,
where $K_1,K_2$ and $K_3$ are finite rank operators.
$K_1$ is the orthogonal projection onto the space of $L^2$-harmonic $p$-forms and comes from the discrete
part of the spectrum. $K_2$ corresponds to the part of the continuous spectrum spanned by co-closed
generalized eigenforms. Using a deformation of the contour of integration in the spectral representation
of the resolvent  one can show (see e.g. Prop. 4.2 in \cite{Str05}) that the integral kernel $k_2(x,y)$ of $K_2$ is given by
$$
 k_2(x,y) = -\frac{\I}{4} \sum_{i} F(\psi_i,0)(x) \overline{F(\psi_i,0)(y)}, 
$$
where $\psi_i$ is an orthonormal basis in $\mathcal{H}^p(Y)$. The part $K_3$ comes from the part of the
continuous spectrum spanned by the closed generalized eigenforms. In the same way as for
co-closed generalized eigenforms one obtains
$$
 k_3(x,y) =-\frac{\I}{4} \sum_{i} F(du \wedge \psi_i,0)(x) \overline{F(du \wedge \psi_i,0)(y)}, 
$$
where $\psi_i$ is an orthonormal basis in $\mathcal{H}^{p-1}(Y)$.
\end{rem}

\section{Cohomology and Hodge theory on $M$} \label{three}

As before let $M$ be a compact manifold with boundary $Y$ and 
$X= (\rz^+ \times Y) \cup_Y M$ the associated  manifold with a cylindrical
end. We consider the long exact sequence 
\begin{equation}\label{exactsequ}
\begin{xy}
\xymatrix{
\ldots\ar[r]^(0.3)\partial & H^k(M,Y,\rz) \ar[r]^(0.5)e & H^k(M,\rz) \ar[r]^{r} 
& H^k(Y,\rz) \ar[r]^(0.4){\partial} & H^{k+1}(M,Y,\rz) \ar[r]^(0.6)e & \ldots,
}
\end{xy}
\end{equation}
in de Rham cohomology. Here $e$ is the canonical embedding and $r$ is the 
restriction homomorphism. There are three cochain complexes which compute the 
relative de Rham  cohomology. Let
\[
\Lambda^p(M,Y):=\left\{\omega\in\Lambda^p(M)\colon i^*\omega=0\right\},
\]
where $i\colon Y\to M$ is the inclusion. Since $d$ commutes with $i^*$, we get 
a complex $\Lambda^*(M,Y)$. Its cohomology is denoted by $H^*(M,Y,\R)$.
 There is an exact sequence of complexes
\[
\begin{xy}
\xymatrix{
0\ar[r] &\Lambda^*(M,Y)\ar[r]^j &  \Lambda^*(M)\ar[r]^{i^*} 
& \Lambda^*(Y)\ar[r] &0,
}
\end{xy}
\]
where $j$ is the inclusion map. It gives rise to the long exact sequence 
(\ref{exactsequ}). The connecting homomorphism $\partial$ is defined as 
follows.
Let $[\phi]\in H^k(Y,\rz)$. Extend $\phi$ to a $k$-form $\omega$ on $M$ such
that $\omega=\phi$ in a neighborhood of the boundary. Then 
\begin{gather}
 \partial[\phi]=[d \omega].
\end{gather}
For the second description consider the cochain complex
$\Lambda_{rel}^*(M,Y)$ of the mapping cone of $i^*$ which is defined by
\[
\Lambda_{rel}^p(M,Y):=\Lambda^p(M)\oplus \Lambda^{p-1}(Y)
\]
with differential $d$ given by
\[
d(\omega,\theta)=(d\omega,i^*\omega-d\theta),\quad \omega\in\Lambda^p(M),\;
\theta\in\Lambda^{p-1}(Y).
\]
Let 
\[
\alpha\colon \Lambda^{p-1}(Y)\to \Lambda^p_{rel}(M,Y)\quad\mathrm{and}\quad
\beta\colon \Lambda^p_{rel}(M,Y)\to \Lambda^p(M)
\]
be defined by $\alpha(\theta)=(0,(-1)^{p-1}\theta)$ and $\beta(\omega,\theta)=
\omega$, respectively. Then $\alpha$ and $\beta$ are cochain maps and we get a
second exact sequence of cochain complexes
\[
\begin{xy}
\xymatrix{
0\ar[r] &\Lambda^{*-1}(Y)\ar[r]^\alpha &  \Lambda^*_{rel}(M,Y)\ar[r]^{\beta} 
& \Lambda^*(M)\ar[r] &0
}
\end{xy}
\]
There is a natural inclusion of cochain complexes
\[
\gamma\colon \Lambda^*(M,Y)\to \Lambda^*_{rel}(M,Y), \quad \omega\mapsto 
(\omega,0).
\]
It follows from the corresponding commutative diagram of long exact sequences
that $\gamma$ induces an isomorphism
\[
\gamma\colon H^*(M,Y,\R)\cong H^*_{rel}(M,Y,\R).
\]
Finally $H^*_{rel}(M,Y,\R)$ is also naturally isomorphic to the cohomology
with compact supports $H_c^*(X)$. The isomorphism can be described as follows.
Let $p\colon Z=\R^+\times Y\to Y$ be the canonical projection. Integration over
the fibre $\R^+$ of $p$ induces a mapping
\[
p_*\colon \Lambda_c^p(\R^+\times Y)\to \Lambda^{p-1}(Y).
\]
Define a map
\[
\xi\colon\omega\in\Lambda_c^p(X)\mapsto (\omega|_M,-p_*(\omega|_Z))\in 
\Lambda^p_{rel}(M,Y)
\]
This is a chain map. If the support of $\omega$ is contained in $M\setminus Y$,
then $\xi(\omega)=(\omega,0)$. Since every cohomology class in $H^p_c(X)$ has
a representative of this form, it follows that $\xi$ induces an isomorphism
\[
\bar\xi\colon H_c^p(X)\to H^p_{rel}(M,Y).
\]

If we fix a metric on $Y$ we may identify $H^k(Y,\rz)$ with the space of 
harmonic forms $\mathcal{H}^k(Y,\rz)$.
In fact, the image of $e$, i.e. the kernel of $r$ can be read off from
the scattering matrix at $0$.

\begin{theorem} \label{imageofS}
 The $+1$ eigenspace of the scattering matrix  $S(0)$ on $H^p(Y,\rz)$ 
coincides with 
 $\mathrm{Im}(H^p(M,\rz) \to H^p(Y,\rz))$.
\end{theorem}
\begin{proof}
 Let $\psi \in \mathcal{H}^p(Y,\rz)$ with $\psi=S(0) \psi$. By 
Proposition \ref{limiting}, $F(\psi,0)$  is closed and coclosed.
Therefore the restriction of $F(\psi,0)$ to $M$  defines a
cohomology class in $H^p(M)$. Expand $F(\psi,0)$ on $Z$ in terms of an 
orthonormal basis of $\ker\Delta_p^\prime\oplus\ker\Delta_{p-1}^\prime$. Using
that  $F(\psi,0)$ is  closed  and coclosed, 
it follows that its expansion on $Z$ has the form
 \begin{gather}\label{expans6}
  F(\psi,0)= 2\psi + \sum_{\mu_{\phi_i}>0} a_i e^{- \mu_{\phi_i} u} 
(d^\prime \phi_i - \mu_{\phi_i} du \wedge \phi_i),
 \end{gather}
 where $\{\phi_i\}_{i\in\nz}$ is an orthonormal basis of $\delta(\Lambda^p(Y))$
consisting of eigenforms of $\Delta^\prime_{p-1}$ with eigenvalues 
$\mu_{\phi_i}^2$. In particular 
 \begin{gather}
   i_Y^*F(\psi,0)= 2\psi + d^\prime\left( \sum_i a_i \phi_i \right),
 \end{gather}
 and therefore the image of the cohomology class $[\frac{1}{2}F(\psi,0)]$ 
under $r$
 is precisely $[\psi]$. Therefore we have shown that the $+1$ eigenspace
 of $S(0)$ is contained in the image of $H^p(M,\rz)$ in $H^p(Y,\rz)$. 
 Now let $\phi$ be an element in the image of $r$, i.e. $\phi$ is a harmonic 
form that is in the same cohomology class as the restriction
 of a closed form $f$ on $M$. If $\psi$ is in the $-1$ eigenspace of $S(0)$
 then by Theorem \ref{main22}, $*\psi$ is in the $+1$ eigenspace and we have
 \begin{gather}
  0=\int_M df \wedge F(*\psi,0) = \int_Y i^*_Y(f) \wedge *\psi = \int_Y \phi 
\wedge *\psi = \langle \phi,\psi \rangle.
 \end{gather}
Hence, any element in $\mathrm{Im}(H^p(M,\rz) \to H^p(Y,\rz))$ is in the
orthogonal complement to the $-1$ eigenspace of $S(0)$
which is exactly the $+1$ eigenspace. This shows that $\mathrm{Im}(H^p(M,\rz) 
\to H^p(Y,\rz))$ is contained in the $+1$ eigenspace and this concludes the 
proof.
\end{proof}

Hence, the scattering matrix at $0$ is determined completely by the metric on 
the boundary.
Namely, it is equal to $1$ on the kernel of $\partial$ and equal to $-1$ on 
the orthogonal
complement of the kernel of $\partial$. 
Recall that by Proposition \ref{limiting}, $F(\psi,0)$ is a closed and 
coclosed $p$-form.
Let  $\hat F: H^k(Y,\rz) \to H^k(M,\rz)$ be the map defined by
\begin{gather}
 \hat F( \psi)= [\frac{1}{2}F(\psi,0)|_M].
\end{gather}
Then, by construction, $r \circ \hat F$ is the orthogonal projection onto the 
kernel of $\partial$.

Hodge theory for manifolds with boundary shows that absolute and
relative cohomology classes have unique harmonic representatives that satisfy
certain boundary conditions. We recall  the definition of the 
relative and absolute boundary conditions
for the Laplace operator. The operator $\Delta_{rel}$ is the closure of the
Laplace operator with respect to the relative boundary conditions
\begin{gather*}
 \omega |_Y=0, \quad (\delta \omega)|_Y=0.
\end{gather*}
The operator $\Delta_{abs}$ is the closure of the
Laplace operator with respect to the absolute boundary conditions
\begin{gather*}
 (*\omega)|_Y=0, \quad (* d \omega)|_Y=0.
\end{gather*}
Both operators are self-adjoint and have compact resolvents. 
Their kernels are the space of harmonic forms satisfying relative and absolute
boundary conditions. Equivalently, they are given by
\begin{gather*}
 \mathcal{H}^p_{rel}(M)=\{\omega \in \Lambda^p(M) \colon d \omega = 
\delta \omega
 =0,\; \omega|_Y=0\},\\
 \mathcal{H}^p_{abs}(M)=\{\omega \in \Lambda^p(M) \colon d \omega = \delta 
\omega=0,\; (*\omega)|_Y=0\},
\end{gather*}
Hodge theory for
manifolds with boundary shows that the canonical maps
\begin{gather*}
 \mathcal{H}^p_{rel}(M) \to H^p(M,Y,\rz),\quad 
 \mathcal{H}^p_{abs}(M) \to H^p(M,\rz)
\end{gather*}
are isomorphisms, that is, every absolute/relative cohomology class has a
unique harmonic representative satisfying absolute/relative boundary
conditions (see e.g. \cite{DS52}).

The harmonic representative $\phi$ of the cohomology class $[\phi] 
\in H^p(M,Y,\rz)$
is the unique minimizer of the functional
\begin{gather*}
 \omega \mapsto \langle \omega, \omega \rangle_{L^2(M)}
\end{gather*}
in $[\phi]$. Similarly, any harmonic form satisfying absolute boundary
conditions minimizes the $L^2$-norm in its absolute cohomology class.
Apart from these minimax principles there is another interesting
minimizing problem which is described in the following proposition.

\begin{theorem} \label{mini}
Let $\phi \in \Lambda^p(Y)$. Consider the functional  $F$ on 
$\{ \omega \in \Lambda^p(M) \colon \omega|_Y=\phi\}$ which is defined by
\[
F(\omega)=\langle d \omega, d \omega \rangle_{L^2}.
\]
Then there exists a unique coclosed harmonic form $\omega_0$ with 
$\omega_0|_Y=\phi$ such that $\omega_0$ is orthogonal to 
$\mathcal{H}^p_{rel}(M)$. The minimum of $F$ is attained at $\omega_0$ and
$\omega_0$ is the unique coclosed minimizer that is orthogonal to 
$\mathcal{H}^p_{rel}(M)$.
If $\phi$ is closed $d \omega_0$ is the harmonic representative in 
$\partial [\phi]$.
\end{theorem}
\begin{proof}
 We divide the proof into several steps.\\
 {\bf Uniqueness:}
 If two forms $\omega_0$ and $\omega_0'$ are harmonic, coclosed, and their
restrictions to $Y$ coincide, it follows that $\omega_0-\omega_0'$
 is harmonic, coclosed and satisfies relative boundary conditions. Therefore,
 $\omega_0-\omega_0' \in \mathcal{H}^p_{rel}(M)$. If both $\omega_0$ and 
$\omega_0'$ are orthogonal to $\mathcal{H}^p_{rel}(M)$ it  follows 
that $\omega_0=\omega_0'$.\\
 {\bf Existence:}
 Choose any extension $\tilde \psi$ of $\phi$ to $M$ which in 
 a neighborhood of $Y$ is of the form
 \begin{gather}
  \phi + u du \wedge \delta \phi.
 \end{gather}
 Then $\delta \tilde \psi$ vanishes near $Y$. 
 Next we claim that the form $\Delta \tilde \psi$ is in the orthogonal 
complement of the kernel
 of $\Delta_{rel}$. Indeed, if $\xi \in \mathcal{H}^p_{rel}(M)$, then
 \begin{gather*}
  \langle \xi, \Delta \tilde \psi \rangle = \int_M d \delta \tilde 
\psi\wedge * \xi
  +\int_M \xi \wedge * \delta d \tilde \psi=\\=
  \int_Y \delta \tilde \psi \wedge *\xi - \int_Y \xi\wedge * d \tilde 
\psi=0,
 \end{gather*}
 where the first integral vanishes because $\tilde \psi$ is coclosed near $Y$
 and the second integral vanishes because $\xi$ satisfies relative boundary
 conditions. Let us denote by $p_0$ the orthogonal projection onto the kernel
 of $\Delta_{rel}$.
 Since $\Delta \tilde \psi$ is in the orthogonal complement of the kernel of 
 $\Delta_{rel}$ the following form is well defined and harmonic
 \begin{gather*}
  \omega_0  = \tilde \psi - 
 \left(\Delta_{rel}|_{\mathcal{H}^p_{rel}(M)^{\perp}} \right)^{-1} 
 (\Delta \tilde \psi) - p_0 \tilde \psi.
 \end{gather*}
 By construction it is also in the orthogonal complement of
 $\mathcal{H}^p_{rel}(M)$ and satisfies $\omega_0|_Y=\phi$.
 Moreover, $\omega_0$ is harmonic and since $\Delta_{rel}$ commutes with
 $\delta$, we have
\[
\delta \omega_0=\delta\tilde \psi - 
 \left(\Delta_{rel}|_{\mathcal{H}^p_{rel}(M)^{\perp}} \right)^{-1} 
 \Delta \delta\tilde \psi 
\]
Since $\delta\tilde \psi$ vanishes near $Y$, it is in the domain of 
$\Delta_{rel}$ and therefore, the right hand side vanishes.
Hence $\delta\omega_0=0$.\\
{\bf Minimizing property:}
 Suppose that $\psi$ is a $p$-form with $\psi|_Y=0$, then
 \begin{gather*}
  \langle d (\omega_0 + \psi),d(\omega_0+ \psi) \rangle = 
   \langle d \omega_0 ,d \omega_0 \rangle + 2 \langle d \omega_0 ,d \psi \rangle
   + \langle d \psi ,d \psi \rangle=\\=
  \langle d \omega_0 ,d \omega_0 \rangle + \langle d \psi ,d \psi \rangle,
 \end{gather*}
 because $\langle d \omega_0 ,d \psi \rangle= \langle \delta d \omega_0 , \psi
 \rangle=0$. Thus, $\omega_0$ is a minimizer. One can see immediately that the
 Euler-Lagrange equations for the minimizer are the equations $\delta d
 \omega_0=0$. Thus, any coclosed minimizer has to be harmonic. The uniqueness
 statement for the minimizer thus follows from the above uniqueness statement.
\end{proof}

Since restriction to the boundary commutes with the differential the map which
sends $\phi$ to $\omega_0$ commutes with the differential. Therefore, if
$\phi$ is exact or closed, so is $\omega_0$.

Now consider the manifold $M_a$ which is obtained from $M$ by attaching the
cylinder $[0,a]\times Y$ to $M$. 
Then $M_a$ is a manifold with boundary $Y_a$. Let  $\phi\in
\mathcal{H}^{p-1}(Y)$. We regard it as a harmonic form on $Y_a$.
By the Hodge theorem there is a unique harmonic form in
$\mathcal{H}^p_{rel}(M_a)$ which represents $\partial [\phi]\in
 H^p(M_a,\partial M_a)$. Let us denote
this form by $\partial_a \phi$, where the notation $\partial_a$ 
indicates that  $\partial_a$ maps $\mathcal{H}^{p-1}(Y)$ to different spaces 
depending on $a$.
For each $a$ we have the $L^2$-inner product on  $\mathcal{H}^p_{rel}(M_a)$. 
It is a natural question to ask whether $\partial$  as a map
from one Hilbert space to another one is a partial isometry. The scalar
product on $\mathcal{H}^p_{rel}(M_a)$, however, depends on $a$. 

\begin{theorem} \label{theoremquadratic}
 Let $Q_a$ be the sesquilinear form on $\mathrm{ker}(\partial)^{\perp}$ which is
defined by
 \begin{gather*}
  Q_a(\psi,\phi)=\langle \partial_a \psi, \partial_a \phi \rangle_{L^2(M_a)}
 \end{gather*}
 and let $q(a)$ be the unique linear operator
 $\mathrm{ker}(\partial)^{\perp} \to \mathrm{ker}(\partial)^{\perp}$ such that
 \begin{gather*}
  Q_a(\psi,\phi)=\langle \psi, q(a) \phi\rangle.
 \end{gather*}
 Then, as $a \to \infty$:
 \begin{gather*}
  q(a)^{-1}= a \id + \frac{\I}{2} S'(0)|_{\mathrm{ker}(\partial)^{\perp}} + 
O( a e^{-\mu_1 a}).
 \end{gather*}
\end{theorem}
\begin{proof}
Recall that $\ker(\partial)^{\perp}=\ker(\mathrm{I}+S(0))$. Let
$\phi\in\ker(\partial)^{\perp}$. Then
by (\ref{expansion2}) the restriction of $F'(\phi,0)$ 
to the cylinder $Z$ has the following form
\begin{gather}
  F'(\phi,0)|_Z= -2 \I u \phi + S'(0) \phi +R'(\phi,0).
 \end{gather}
It follows that $F'(\phi,0)$ is a coclosed harmonic form. 
 Thus, by theorem \ref{mini} $F'(\phi,0)$ is a coclosed minimizer
 of the functional $\eta \to \langle d\eta,d \eta\rangle_{L^2(M_a)}$
 with boundary condition $\eta|_{Y_a}=-2 \I a \phi + S'(0) \phi +
 R'(\phi,0)|_{Y_a}$. Let $H$ be the minimizer with boundary conditions
 $H|_Y=a \phi + \frac{\I}{2} S'(0)\phi$ so that $d H$ is the unique harmonic
 representative of  $\partial_a \left(( a \id + \frac{\I}{2} S'(0))\phi\right)$.
 Again by theorem \ref{mini} 
 \begin{gather}
   G_\phi:=\frac{i}{2}F'(\phi,0)-H
 \end{gather}
minimizes the functional $\eta \to \langle d\eta,d \eta\rangle_{L^2(M_a)}$
with boundary conditions
\begin{gather}
  \eta|_{Y_a}=R'(\phi,0)|_{Y_a}.
\end{gather}
Then for every $\psi\in\ker(\partial)^{\perp}$ we have
\begin{gather} \label{eqncxcx}
 Q_a(\psi,( a \id + \frac{\I}{2} S'(0))\phi )= 
\langle \partial_a \psi,\frac{i}{2}d F'(\phi,0)
 \rangle_{L^2(M_a)} - \langle \partial_a \psi,d G_\phi \rangle_{L^2(M_a)}. 
\end{gather}
 The second term on the right hand side can be estimated using the
 Cauchy-Schwarz inequality
 \begin{gather}
  |\langle \partial_a \psi,d G_\phi \rangle_{L^2(M_a)}| \leq \|\partial_a \psi\|_{L^2(M_a)}
 \cdot \|dG_\phi \|_{L^2(M_a)}.
\end{gather}
The first term in (\ref{eqncxcx}) can be explicitly calculated. Using 
(\ref{fprime}) we get
\begin{gather*}
 \langle \partial_a \psi,\frac{i}{2}d F'(\phi,0) \rangle_{L^2(M_a)}=
\frac{1}{2}\int_Y \psi \wedge *_M F(du \wedge \phi,0)=\\
  =\int_Y \psi \wedge * \phi=\langle \psi, \phi  \rangle,
\end{gather*}
where we used that $R(du \wedge \phi,0)|_Y$ is orthogonal to $\psi$.
Thus,
\begin{gather} \label{eqnghh}
 |Q_a(\psi,( a \id + \frac{\I}{2} S'(0))\phi )- 
\langle\psi,\phi\rangle| \leq
 \|\partial_a \psi\|_{L^2(M_a)} \cdot \|dG_\phi\|_{L^2(M_a)}.
\end{gather}
The terms on the right hand side can be estimated as follows. First note that
 $\|dG_\phi\|_{L^2(M_a)}$ minimizes $\|d\eta\|_{L^2(M_a)}$ over all forms
 $\eta$ which restrict to $R'(\phi,0)|_{Y_a}$ on $Y_a$. Moreover,
 $\chi_a:=R'(\phi,0)|_{Y_a}$ is exponentially decaying in $a$.
 Define the form $\eta_a$ by $\eta_a:=\frac{u}{a}\chi_a$
 on the cylinder and $0$ elsewhere. Then,
 \begin{gather}
  \| d \eta_a \|_{L^2(M_a)}^2 = \frac{a}{3} \|d \chi_a\|^2_{L^2(Y)} +
  \frac{1}{a} \|\chi_a\|^2_{L^2(Y)}.
 \end{gather}
By Lemma \ref{rapiddec} we have
 \begin{gather}
  \|d \chi_a\|^2_{L^2(Y)} + \|\chi_a\|^2_{L^2(Y)} \leq C_\phi e^{- 2\mu_1 a},
 \end{gather}
which implies
 \begin{gather}
  \|dG_\phi\|_{L^2(M_a)} \leq \tilde C_\phi e^{- \mu_1 a}.
 \end{gather}
To estimate $\parallel\partial_a\psi\parallel_{L^2(M_a)}$, recall that by
Theorem \ref{mini}, 
$\partial_a\psi=d\omega_0$, where $\omega_0$ 
minimizes the functional $\eta\mapsto\| d\eta\|_{L^2(M_a)}^2$ with boundary 
conditions $\eta|_{Y_a}=\psi$. Let $f\in C^\infty(\rz^-)$ such that $f(u)=1$
for $-1/4\le u\le 0$ and $f(u)=0$ for $u\le -3/4$. For $a\ge 1$ define
$\hat\psi_a\in \Lambda^p(M_a)$ by
\[
\hat\psi_a(x)=\begin{cases} f(u-a)\psi(y),& \mathrm{if}\;
 x=(u,y)\in [0,a]\times Y;\\
0,& \mathrm{otherwise}.
\end{cases}
\]
Then it follows that there exists $C>0$ such that
\[
\|\partial_a\psi\|_{L^2\Lambda^{p+1}(M_a)}\le \|d\hat\psi_a\|_{L^2\Lambda^{p+1}(M_a)}
\le C\|\psi\|_{L^2\Lambda^p(Y)}.
\]
Thus we get
\begin{gather} \label{eqngfh}
 |Q_a(\psi,( a \id + \frac{\I}{2} S'(0))\phi )- \langle\psi,\phi\rangle| \leq
 C_\phi'' \|\psi\|e^{- \mu_1 a}.
\end{gather}
By compactness of the unit sphere in $\mathrm{ker}(\partial)^{\perp}$
the constant can be chosen independent of $\phi$ if we use vectors of norm 
$1$ only. Hence
\begin{gather*}
 \sup_{\|\psi\|,\|\phi\|=1} \left\langle \psi , (q(a)(a \id + 
\frac{\I}{2} S'(0)) - \id)\phi \right\rangle \leq C e^{- \mu_1 a},
\end{gather*}
which implies
\begin{equation}\label{eqnfin}
\|a \id + \frac{\I}{2} S'(0)-q(a)^{-1}\|\le C\|q(a)^{-1}\|e^{- \mu_1 a}.
\end{equation}
From this inequality we deduce that $\|q(a)^{-1}\|\le C(1+a)$. Combined with
(\ref{eqnfin}) the statement of the theorem follows.
\end{proof}

\section{Hodge theory on $X$ and the scattering length} \label{four}

In this section we give a description of the long exact cohomology sequence of
$X$ in terms of harmonic forms and we  derive a cohomological formula
for the scattering length.

Let $H_c^*(X)$ denote the de Rham cohomology groups with compact supports.
It is well known (see \cite{Mel93}, Sec. 6.4) that $H^*(X)$ and $H^*_c(X)$ 
are canonically isomorphic to certain 
spaces of extended harmonic forms on $X$. We recall some details. 

The space of extended harmonic forms $\mathcal{H}_{ext}^p(X)$ is 
defined to be the subspace of all (real valued) $\psi \in \Lambda^p(X)$ satisfying 1) 
$\Delta_p \psi=0$ 
and 2) there  exist $\phi_1 \in \ker \Delta_p'$ and $\phi_2 \in \ker 
\Delta_{p-1}'$ such that 
\begin{gather*}
\psi|_Z - \phi_1 - du \wedge \phi_2 \in L^2\Lambda^p(Z).
\end{gather*}

Note that for a  given $\psi \in \mathcal{H}_{ext}^p(X)$ the sections
$\phi_1$ and $\phi_2$ are uniquely determined. We regard $\phi_1$ 
(resp. $\phi_2$)  as the tangential (resp. normal) boundary value of $\psi$ 
at infinity and we denote them by $\psi_t$ and $\psi_n$, respectively.
The spaces satisfying absolute and relative boundary
conditions at infinity are then defined as
\begin{gather*}
 \mathcal{H}_{ext,abs}^p(X):=\{\psi \in \mathcal{H}_{ext}^p(X) \mid
 \psi_n=0\},\\
\mathcal{H}_{ext,rel}^p(X):=\{\psi \in \mathcal{H}_{ext}^p(X) \mid
 \psi_t=0\}.
\end{gather*}
Since   $\psi\in\mathcal{H}_{ext}^p(X)$ is harmonic and 
the form $\psi - \psi_t - du \wedge \psi_n$ is square integrable, it follows 
from (\ref{expansion1}) that there exists $c>0$ such that
\begin{equation}\label{estimation8}
(\psi - \psi_t - du \wedge \psi_n)(u,y)\ll e^{-cu},\quad (u,y)\in Z.
\end{equation}
Moreover  $d\psi$ and $\delta\psi$ are also exponentially decaying. Applying
Greens formula to $M_a$, we get 
\[
0=\langle\Delta\psi,\psi\rangle_{M_a} =\parallel d\psi\parallel_{M_a}^2+
\parallel\delta\psi\parallel_{M_a}^2+O(e^{-ca}),
\]
which implies that 
\begin{equation}\label{closed}
d\psi=0,\;\delta\psi=0\quad\mathrm{for}\;\mathrm{all}\;\psi\in\cH^p_{ext}(X).
\end{equation}
The intersection $\mathcal{H}_{ext,abs}^p(X) \cap \mathcal{H}_{ext,rel}^p(X)$
is the space $\mathcal{H}_{(2)}^p(X)$ of square integrable harmonic forms.

On $\mathcal{H}_{ext}^p(X)$ we introduce an inner product as follows. 
For  
$\psi,\phi \in \mathcal{H}_{ext}^p(X)$ let
\begin{equation}\label{innprod1}
\langle \psi, \phi \rangle = \int_M \psi \wedge * \phi +
\int_Z(\psi-\psi_t-du\wedge\psi_n)\wedge*{
(\phi-\phi_t-du\wedge \phi_n)}.
\end{equation} 
To verify that this is an inner product, we only need to show that 
$\parallel\phi\parallel=0$ implies $\phi=0$. So suppose that $\parallel\phi
\parallel=0$. Then, in particular, we have $\phi|_M=0$ and the unique
continuation property for harmonic forms (see e.g. \cite{BB00,BW93})
implies $\phi=0$. We note that the 
inner product can be also defined by the following formula:
\begin{equation}\label{innprod2}
\langle \psi, \phi \rangle =\lim_{a\to\infty}\left(\int_{M_a}\psi\wedge *\phi
-a(\langle\psi_t,\phi_t\rangle+\langle\psi_n,\phi_n\rangle)\right).
\end{equation}

This inner product coincides on the subspace $\cH^p_{(2)}(X)$ with the usual 
inner product on $\cH^p_{(2)}(X)$. The orthogonal projections define canonical
maps
\begin{equation*}
 \cH_{ext,rel}^p(X) \to \cH_{(2)}^p(X),\quad
 \cH_{ext,abs}^p(X) \to \cH_{(2)}^p(X).
\end{equation*}
Moreover, we have the maps
\begin{equation}\label{maps1}
\begin{split}
 & \hat F: \mathcal{H}^p(Y) \to  \mathcal{H}_{ext,abs}^p(X), \quad \phi \mapsto
  \frac{1}{2} F(\phi,0),\\
 &\hat G: \mathcal{H}^p(Y) \to  \mathcal{H}_{ext,rel}^p(X), 
\quad \phi \mapsto \frac{\I}{2}dF'(\phi,0).
\end{split}
\end{equation}
Next we define maps into the de Rham cohomology. Let 
$\phi\in\cH^p_{ext,abs}(X)$. By (\ref{closed}), $\phi$ is closed and we get a 
canonical map
\begin{equation}\label{derham}
R\colon \cH^p_{ext,abs}(X)\to H^p(X,\R).
\end{equation}
Now consider $\psi \in \mathcal{H}_{ext,rel}^p(X)$. By (\ref{closed}) $\psi$
is closed and on the cylinder $\psi$ is of the form
\begin{equation}\label{restrict}
 \psi|_Z = du \wedge \psi_n + d\theta,
\end{equation}
where $\theta$ is exponentially decaying. 
Let $\chi$ be a function with support on the cylinder $Z$
which is equal to 1 outside a compact set. Following \cite{Mel93}
we can then define a map
\begin{gather}\label{derham1}
R_c\colon \mathcal{H}_{ext,rel}^p(X) \to H^p_c(X,\rz), \quad 
 \psi \mapsto [\psi - d(\chi(u \psi_n + \theta))].
\end{gather}
This map is well defined and independent of the choice of $\chi$.
Indeed changing $\chi$ on a compact subset changes
$$
 \psi - d(\chi(u \psi_n + \theta))
$$
by the differential of a compactly supported form. Let
\[
(\cdot,\cdot)\colon H^p_c(X) \times H^{n-p}(X) \to \rz
\]
be the canonical pairing defined by
\[
([\phi],[\psi])=\int_X\phi\wedge \psi,\quad [\phi]\in H^p_c(X),\; [\psi]\in  
H^{n-p}(X).
\]
Define a pairing 
\begin{gather}\label{pairing}
 (\cdot,\cdot)_{ext}\colon\mathcal{H}_{ext,rel}^p(X) 
\times\mathcal{H}^{n-p}_{ext,abs}(X) \to \rz
\end{gather}
by taking the
constant term in the asymptotic expansion of
\begin{gather*}
 \int_{M_a} \psi \wedge \phi
\end{gather*}
as $a \to \infty$. Applying Green's formula to $M_a$, it follows that there 
exists $c>0$ such that
\begin{gather*}
 \int_{M_a} \psi \wedge  \phi = a \int_{Y} \psi_n \wedge \phi_t +
([\psi - d(\chi(u \cdot \psi_n + \theta))],[\phi] ) + O(e^{-c a}) 
\end{gather*}
as $a \to \infty$. This implies that  the following diagram commutes
{\footnotesize
\begin{equation}\label{diagram}
\begindc{0}[70]
  \obj(20,40)[A]{$\mathcal{H}_{ext,rel}^p(X)$}
  \obj(25,40)[B]{$\times$}
  \obj(30,40)[D]{$\mathcal{H}^{n-p}_{ext,abs}(X)$}
  \obj(40,35)[E]{$\rz$}
  \obj(20,30)[F]{$H^p_c(X,\rz)$}
  \obj(25,30)[G]{$\times$}
  \obj(30,30)[H]{$H^{n-p}(X,\rz)$}
  \mor{D}{E}{}
  \mor{H}{E}{}
  \mor{A}{F}{$R_c$}
  \mor{D}{H}{$R$}
 \enddc,
\end{equation}
}
where the horizontal maps are given by the corresponding pairing.

Let $\phi\in\cH_{ext}^p(X)$ and $\omega\in \Lambda_c^{p-1}(X)$. Since $\phi$
is co-closed (\ref{closed}), it follows that $\langle d\omega,\phi\rangle=0$.
Therefore, for $\phi\in\cH_{ext}^p(X)$, $[\psi]\in H^p_c(X)$, and $\psi^\prime\in
[\psi]$, 
the inner product $\langle \psi^\prime,\phi\rangle$ is independent of the 
representative of the cohomology class $[\psi]$ and will be denoted by
$\langle [\psi],\phi\rangle$. This leads to the following alternative 
description of the inner product in $\cH^p_{ext,rel}(X)$. 

\begin{lem}\label{innpro3}
 For all $\phi,\psi\in\cH^p_{ext,rel}(X)$ we have
\[
\langle \psi,\phi\rangle=\langle R_c(\psi),\phi\rangle.
\]
\end{lem}
\begin{proof}
Applying Stoke's theorem and using that $\theta$ is rapidly decreasing,  we get
\[
\int_{M_a}d(\chi(u\psi_n+\theta))\wedge *\phi=a\langle\psi_n,\phi_n\rangle
+O(e^{-ca}).
\]
By (\ref{innprod2})  we get
\[
\langle \psi,\phi\rangle=\langle \psi-d(\chi(u\psi_n+\theta)),\phi\rangle
=\langle R_c(\psi),\phi\rangle.
\]
\end{proof}

Our next goal is to describe the connecting
homomorphism $\partial: H^p(Y,\rz) \to H^{p+1}_c(X,\R)$ on the level of
harmonic forms. To this end we need some preparation. 
Let $\psi$ be in the $-1$ eigenspace of $S(0)$. Then by (\ref{expansion2}), we
have on $Z$ 
$$
 \frac{\I}{2} F'(\psi,0)|_Z=  u \psi + \frac{\I}{2} S'(0) \psi + \theta,
$$
where $\theta$ is exponentially decaying. Let $\chi$ be a smooth function with
support in $Z$ which is equal to 1 outside a compact set. Then  
$$
 \frac{\I}{2} F'(\psi,0)-\chi( u \psi  +\theta)
$$
is equal to $\frac{\I}{2} S'(0)\psi$ outside a compact set and we conclude that
$$
 d(\frac{\I}{2} F'(\psi,0)-\chi( u \psi  +\theta))
$$
represents $\partial [\frac{\I}{2} S'(0)\psi]$ in $H^{p+1}_c(X,\R)$. 
Let $\kappa\colon \cH^p(Y)\to H^p(Y)$ be the canonical isomorphism. Then we
have shown that for each $\psi\in\cH^p(Y)$ we have
\begin{equation}\label{equ5}
R_c(\hat G(\psi))=\partial\left[\kappa(\frac{i}{2}S^\prime(0)\psi)\right],
\end{equation}
where $\hat G(\psi)$ is defined by (\ref{maps1}). 
\begin{lem}\label{invert} The operator $S^\prime(0)$ in $\cH^*(Y)$ is 
invertible.
\end{lem}
\begin{proof}
Differentiating equations (ii) and (iii) of Theorem \ref{main22}, it
follows that  $S^\prime(0)$ commutes with $S(0)$ and anti-commutes with $*$. 
Therefore, it suffices to show that the restriction of $S^\prime(0)$ to the
$-1$-eigenspace $E_-$ of $S(0)$ is invertible. Let $\psi\in E_-$. Then 
$S^\prime(0)\psi\in E_-$. By Theorem \ref{imageofS} we have 
$E_-=(\ker\partial)^\perp$. Using (\ref{equ5}), it follows that it suffices
to show that $R_c(\hat G(\psi))\not=0$ whenever $\psi\not=0$. By Lemma 
\ref{innpro3} we have
\[
\langle\hat G(\psi),\hat G(\psi)\rangle=\langle R_c(\hat G(\psi)),\hat G(\psi)
\rangle
\]
for all $\psi\in\cH^p_{ext,rel}(X)$. Recall that $\hat G(\psi)$ is a harmonic
form, which is non-zero, if $\psi\not=0$. Therefore, the left hand side of the
above equality is non-zero, if $\psi\not=0$.
\end{proof}

Now we define  maps
\[
\tilde e\colon \cH^p_{ext,rel}(X)\to \cH_{ext,abs}^p(X),\quad \tilde r
\colon \cH_{ext,abs}^p(X)
\to \cH^p(Y),\quad \tilde\partial \colon\cH^p(Y)\to \cH^{p+1}_{ext,rel}(X)
\]
as follows. 
Let $\tilde e$  be the composition of
the orthogonal projection $\cH^p_{ext,rel}(X)\to\cH^p_{(2)}(X)$ and the 
inclusion $\cH^p_{(2)}(X)\to \cH^p_{ext,abs}(X)$. $\tilde r$ 
assigns to $\phi\in\cH_{ext,abs}^p(X)$ its limiting value  $\phi_t$. To define
$\tilde\partial$, we note that by  Lemma \ref{invert}, $S'(0)$ is an
invertible operator. Put
\[
\tilde\partial =\hat G\circ \left( \frac{\I}{2}S'(0) \right)^{-1}.
\]
\begin{pro}
The sequence
{\footnotesize
\[
\begindc{0}[70]
  \obj( 0,20)[D]{$\ldots$}
  \obj(10,20)[A]{$\cH^p_{ext,rel}(X)$}
  \obj(20,20)[B]{$\cH^p_{ext,abs}(X)$}
  \obj(30,20)[C]{$\cH^p(Y)$}
  \obj(40,20)[E]{$\cH^{p+1}_{ext,rel}(X)$}
  \obj(50,20)[F]{$\ldots$}
  \mor{D}{A}{$\tilde\partial$}
  \mor{A}{B}{$\tilde e$}
  \mor{B}{C}{$\tilde r$}
  \mor{C}{E}{$\tilde\partial$}
  \mor{E}{F}{$\tilde e$}
\enddc
\]
}
is exact.
\end{pro}
\begin{proof}
Let $E_\pm=\ker(S(0)\mp\Id)$. By Proposition \ref{limiting} and (\ref{niceform})
it follows that 
\begin{equation}\label{imageofr}
\Im(\tilde r)=E_+=\ker(\hat G).
\end{equation}
 Since $S'(0)$ preserves 
$E_\pm$, we get $\Im(\tilde r)=\ker(\tilde\partial)$. By definition we have
$\Im(\tilde e)=\cH^p_{(2)}(X)$ and this is also equal to $\ker(\tilde r)$. Finally
by Corollary \ref{orthogonal} it follows that $\Im(\tilde\partial)=
\cH^{p+1}_{(2)}(X)^\perp$.
On the other hand, by definition we have $\ker(\tilde e)=
\cH^{p+1}_{(2)}(X)^\perp$.
Thus $\Im(\tilde\partial)=\ker(\tilde e)$.
\end{proof}

 Using the definition of
$\tilde\partial$ and (\ref{equ5}), it follows that
\[
R_c\circ\tilde\partial=\partial\circ \kappa.
\]

By \cite{APS75} every element in
the image of $H^*_c(X,\rz)$ in $H^*(X,\rz)$ can be represented by
a unique square integrable harmonic form. Using these facts we obtain the
 following commutative diagram.

{\footnotesize
$$
\begindc{0}[70]
  \obj( 0,30)[A]{$\ldots$}
  \obj(10,30)[B]{$\mathcal{H}^p_{ext,rel}(X)$}
  \obj(10,40)[CC]{$\mathcal{H}^p(Y)$}
  \obj(20,40)[C]{$\mathcal{H}^p_{(2)}(X)$}
  \obj(20,30)[D]{$\mathcal{H}^p_{ext,abs}(X)$}
  \obj(30,30)[E]{$\mathcal{H}^p(Y)$}
  \obj(30,40)[EE]{$\mathrm{ker}(S(0)-\id)$}
  \obj(40,40)[CCC]{$\mathcal{H}^{p}(Y)$}
  \obj(40,30)[F]{$\mathcal{H}^{p+1}_{ext,rel}(X)$}
  \obj(50,30)[G]{$\ldots$}
  \obj( 0,20)[H]{$\ldots$}
  \obj(10,20)[I]{$H^p_c(X,\rz)$}
  \obj(20,20)[J]{$H^p(X,\rz)$}
  \obj(30,20)[K]{$H^p(Y,\rz)$}
   \obj(40,20)[L]{$H^{p+1}_c(X,\rz)$}
  \obj(50,20)[M]{$\ldots$}
  \mor{H}{I}{$\partial$}
  \mor{I}{J}{$e$}
  \mor{J}{K}{$r$}
  \mor{K}{L}{$\partial$}
  \mor{L}{M}{$e$}
  \mor{A}{B}{$\tilde\partial$}
  \mor{C}{B}{}
  \mor{C}{D}{}
  \mor{D}{E}{$\tilde r$}
  \mor{E}{F}{$\tilde\partial$}
  \mor{F}{G}{$\tilde e$}
  \mor{B}{I}{$R_c$}
  \mor{D}{J}{$R$}
  \mor{E}{K}{$\kappa$}
  \mor{F}{L}{$R_c$}
  \mor{EE}{E}{}
  \mor{EE}{D}{$\hat F$}
  \mor{B}{D}{$\tilde e$}
  \mor{CC}{A}{$\frac{\I}{2} S'(0)$}
  \mor{CC}{B}{$\hat G$}
  \mor{CCC}{E}{$\frac{\I}{2} S'(0)$}
  \mor{CCC}{F}{$\hat G$}
 \enddc
$$
}
\begin{pro}
The maps 
\[R\colon\mathcal{H}_{ext,abs}^p(X) \to H^p(X,\rz),\quad
 R_c\colon \mathcal{H}_{ext,rel}^p(X) \to H^p_c(X,\rz)
\]
are isomorphisms.
\end{pro} 
\begin{proof}
We first consider $R$.
Let $H^p_{!}(X,R)=\Im(e)$. By \cite{APS75}, $R$ induces an isomorphism of 
$\cH^p_{(2)}(X)$ onto $H_!^p(X)$. Let $\phi\in\cH_{ext,abs}^p(X)$ and suppose
that $R(\phi)=0$. Then it follows that $\tilde r(\phi)=0$. Hence $\phi\in
\cH^p_{(2)}(X)$. Since $R$ is an isomorphism on $\phi\in\cH^p_{(2)}(X)$, we get
$\phi=0$. This proves injectivity. Let $\psi\in H^p(X,\R)$. Using 
$H^p(X,\R)\cong H^p(M,\R)$ and Theorem \ref{imageofS}, it follows that 
$\kappa^{-1}(r(\psi))\in \ker(S(0)-\Id)$. Thus by (\ref{imageofr}) there exists 
$\phi\in\cH^p_{ext,abs}(X)$ such that $\tilde r(\phi)=\kappa^{-1}(r(\psi))$. Then
$r(R(\phi)-\psi)=0$. Hence $R(\phi)-\psi\in H_{!}^p(X,\R)$. By the above remark
there is $\omega\in\mathcal{H}^p_{(2)}(X)$  such that $R(\omega)=R(\phi)-\psi$.
Thus $R$
is surjective and hence an isomorphism. Applying the commutativity of the
diagramm and the 5-Lemma, it follows that $R_c$ is an isomorphism too. 
\end{proof}
This can also be proved by
slightly different methods (e.g. \cite{Mel95} and  \cite{Mel93}, Sec.6.4).
\begin{cor}
 In every class in $H^p(X,\rz)$ there is a unique representative in
 $\mathcal{H}^p_{ext,abs}(X)$.
\end{cor}

\begin{cor}\label{extrel}
 For every class $[\psi]$ in $H^p_c(X,\rz)$ there is a unique element 
$\hat \psi$ in
 $\mathcal{H}_{ext,rel}^p(X)$ such that for any $\phi \in
 \mathcal{H}^p_{ext,abs}(X)$:
 \begin{gather*}
  \langle [\psi],[\phi] \rangle = \langle \hat\psi, \phi \rangle.
 \end{gather*}
 Moreover, the map
 \begin{gather*}
  H^p_c(X,\rz) \to \mathcal{H}_{ext,rel}^p(X),\\ [\psi] \mapsto \hat \psi
 \end{gather*}
 is an isomorphism.
\end{cor}

We can now consider the scattering length
\begin{equation}\label{scattleng}
T(0):=-\I S(0)^* S'(0)=-\I S(0) S'(0).
\end{equation}
Let $\partial\colon H^p(Y,\R)\to
H^{p+1}_c(X,\R)$ be the connecting homomorphism. We identify $H^p(Y,\R)$ with
$\cH^p(Y)$ and $H^{p+1}_c(X,\R)$ with $\cH^{p+1}_{ext,rel}(X)$ via Corollary
\ref{extrel}. Thus we may regard the connecting homomorphism as a map
\[
\partial\colon \cH^p(Y)\to \cH^{p+1}_{ext,rel}(X).
\]
Let $(\ker\partial)^\perp$ be
the orthogonal complement of $\ker\partial$. 
\begin{theorem} \label{main3}
The scattering length $T(0)$ is a positive, invertible operator
in $\cH^p(Y)$. It is uniquely determined by the following conditions.
 \begin{gather}\label{scattleng1}
  \forall \phi,\psi \in (\ker{\partial})^{\perp}: 
  \langle \partial \phi , \partial  ( T(0)  \psi) \rangle = 2 \langle \phi,
  \psi \rangle,\\
  T(0) * = * T(0).
 \end{gather}
\end{theorem}
\begin{proof} By Theorem \ref{imageofS}, 
 $\Im(r)=\ker{\partial}$ equals the $+1$ eigenspace of $S(0)$. Therefore
 $(\ker{\partial})^\perp$ equals the $-1$ eigenspace of $S(0)$.
 By Theorem \ref{main22}, $*$ anti-commutes with $S(0)$. Hence
 it interchanges the $\pm1$-eigenspaces.
 $*$ also anti-commutes with $S'(0)$ and consequently
 the scattering length $T(0)$ commutes with the Hodge star operator.
 It is therefore completely determined by its restriction to
 $(\ker{\partial})^\perp$. Let $\phi,\psi\in (\ker{\partial})^\perp$. Using
(\ref{equ5}) and Lemma \ref{innpro3}, we have
\begin{equation}
\begin{split}
\langle\partial\phi,\partial T(0)\rangle&=\langle\partial\phi,-i\partial
(S^\prime(0)S(0)\psi)\rangle=2\langle\partial\phi,\partial
\left((i/2)S^\prime(0)\psi\right)\rangle\\
&=2\langle\partial\phi,R_c(\hat G(\psi))\rangle=2\langle\partial\phi,
\hat G(\psi)\rangle.
\end{split}
\end{equation}

Let 
$\tilde\phi$ be the pull-back of $\phi$ to $Z$ and let $\chi\in C^\infty(Z)$
such that $\chi=0$ on $[0,1]\times Y$ and $\chi=1$ outside a compact set.
Then $d(\chi\tilde\phi)\in\Lambda^{p+1}_c(X)$ represents $\partial[\phi]\in
H^{p+1}_c(X,\R)$. By the same argument as in the proof of Lemma \ref{innpro3} 
we get
\[
\langle \partial \phi, \hat G(\psi) \rangle=
\langle d(\chi\tilde\phi),\hat G(\psi)\rangle.
\]
By (\ref{limitvalue}), 
$*\hat G(\psi)=*\frac{i}{2}dF'(\psi,0)$ is an extended harmonic form with 
limiting value $*\psi$. Using Stokes theorem, applied to $M_a$, it follows that
\begin{equation}\label{innerpr3}
\langle d(\chi\tilde\phi),\hat G(\psi)\rangle=
\lim_{a\to\infty}\int_{\partial M_a}\phi\wedge *\hat G(\psi)
=  \int_Y \phi \wedge * \psi =\langle \phi, \psi \rangle.
\end{equation}
This concludes the proof of the theorem.
\end{proof}
\begin{cor}
 For the scattering length $T(0)=-\I S(0)^* S'(0)$ we have the following 
formula
\[
T^{-1}(0) = \frac{1}{2}\left(\partial^* \partial + (*)^{-1}\partial^* \partial
   * \right)=\frac{1}{2}\left( r r^* + (*)^{-1} r r^* * \right).
\]
\end{cor}
This implies that $\frac{1}{2} r^* T(0) r$ is equal to the orthogonal
projection onto the orthogonal complement of $\ker r$.

\section{Estimates on the norm of extended harmonic forms} \label{lzweiestimates}

By the results of the previous sections we have canonical isomorphisms
\begin{gather}\label{caniso}
\begin{split}
&\eta_{abs}\colon  \mathcal{H}^p_{ext,abs}(X) \cong H^p(X,\R)\cong H^p(M,\R),\\
&\eta_{rel}\colon  \mathcal{H}^p_{ext,rel}(X) \cong H^p_c(X,\R)\cong H^p(M,Y,\R).
\end{split}
\end{gather}
In this section we establish relations between some norms on
these spaces. On the cohomology groups $H^p(M,\rz)$ and $H^p(M,Y,\rz)$
there is the so-called comass norm (see \cite[Ch. 4C]{Gro99}) which is defined as follows.
If $V$ is a finite dimensional inner product space, then $\Lambda^p V^*$ has a natural inner
product as well and we denote the norm that is induced by this inner product
by $\|\cdot\|$. The comass norm $\| \cdot \|_\infty$ on $\Lambda^p V^*$ is defined by
\begin{gather}
 \| \omega \|_\infty = \mathrm{sup} \{\omega(e_1,\ldots,e_p) \mid \: e_k \in
 V, \|e_k\|=1\}
\end{gather}
Since the norms are equivalent there is a constant $C$ such that
\begin{gather} \label{estconst}
 \| \omega \|^2 \leq C \| \omega \|_\infty^2,
\end{gather}
and we denote by $C(n,p)$ the optimal such constant.
Since all $n$-dimensional inner product spaces are unitarily equivalent the
constant depends only on $n$ and $p$. Of course (see also \cite{Fed69}),
\begin{gather}
 C(n,0)=C(n,1)=1,\\
 C(n,p) \leq {n \choose p}.
\end{gather}
Moreover, since the Hodge star operator leaves the space of primitive
forms invariant, we have
\begin{gather}
 C(n,n-p)=C(n,p).
\end{gather}
It is also known that
\begin{gather}
 C(n,2)=[\frac{n}{2}]
\end{gather}
(see \cite{GlKo02}). Now let $B$ be a differentiable manifold. Let $\omega\in
\Lambda^p(B)$. 
The comass $\| \omega \|_\infty$ of  $\omega$ is defined by
\begin{gather}
 \| \omega \|_\infty =
 \sup \{ \omega_x(e_1,\ldots,e_p)  \mid \: x \in B, \: e_i \in T_x B,\:
 g(e_i,e_i)=1 \}=\\=\sup \{ \|\omega_x\|_\infty \mid \: x \in B \}. \nonumber
\end{gather}
For a compact manifold $B$ with smooth boundary $\partial B$ 
 this induces a norm on $H^p(B,\partial B,\rz)$ by
\begin{gather}
 \| \phi \|_\infty = \inf \{ \|\omega\|_{\infty} \mid \phi = [\omega],\;
 \omega\in\Lambda^p(B,\partial B),\,d\omega=0\},
\end{gather} 
which we also refer to as the comass norm. To compare the norms on the various
cohomology groups, we need some preparation.
Let $\psi \in \Lambda^p(M)$. 
We define an extension $\hat \psi \in L^\infty\Lambda^p(X)$ of $\psi$
in the following way. The restriction $\psi|_Y$ can be expanded into
eigensections of the Laplace-Beltrami operator on $Y$:
\begin{gather*}
 \psi|_Y= \phi + \sum_{i=1}^\infty a_i \phi_i,
\end{gather*}
where $\phi$ is harmonic and $\phi_i$ is an orthonormal basis in the orthogonal complement
of the space of harmonic forms such that
\begin{gather*}
 \Delta' \phi_i = \mu_{\phi_i}^2 \phi_i.
\end{gather*}
Now define
\begin{gather}\label{extension1}
 \hat \psi(x) = \left \{ \begin{matrix} \psi(x) \quad \textrm{for } x \in M,\\ 
 \phi(y) + \sum_{i=1}^\infty a_i e^{-\mu_{\phi_i} u } \left( \phi_i(y) - \mu_{\phi_i}^{-1} du \wedge
 \delta' \phi_i(y) \right) \quad \textrm{for } x=(u,y) \in Z \end{matrix} \right.
\end{gather}
As usually $x=(u,y)$.
The map $\psi \mapsto \hat \psi$ is of course linear and maps
into the space of bounded sections. Note that in general $\hat \psi$
is not continuous. However, it satisfies 
\[
i_Y^*(\hat\psi|_M)=i^*_Y(\hat\psi|_Z).
\]
Therefore, using Green's formula, it follows that the distributional 
derivative $d \hat \psi$ satisfies
\begin{gather}\label{derivative}
 d \hat \psi(x) = \sum_{i=1}^\infty a_i e^{-\mu_{\phi_i} u} \left( d' \phi_i 
- \mu_{\phi_i}
   du \wedge \phi_i + \mu_{\phi_i}^{-1} du \wedge d' \delta' \phi_i \right)=\\
 \nonumber =
 \sum_{i=1}^\infty a_i e^{-\mu_{\phi_i} u} \left( d' \phi_i - \mu_{\phi_i}^{-1} du 
\wedge\delta' d' \phi_i \right)= \widehat{d \psi}.
\end{gather}
In particular $d \hat \psi$ is again bounded. If $\psi$ is closed, then
$\hat \psi$ is also closed. Note that the extension map 
$ \psi\in\Lambda^p(M)\mapsto
\hat\psi\in L^\infty\Lambda^p(X)$ is chosen so that it inverts the restriction
operator on the space $\mathcal{H}^p_{ext,abs}(X)$. Namely, for $F \in
\mathcal{H}^p_{ext,abs}(X)$ it follows from (\ref{expans6}) that
\begin{gather}
 F=\widehat{F|_M}.
\end{gather}
We can now establish the comparison results for the norms.
\begin{lem} \label{hilfreichlemma}
 Let $\psi \in \Lambda^p(M)$ be closed 
 and $F \in \mathcal{H}^p_{ext,abs}(X)$ a harmonic form such that $F|_M$ and 
$\psi$ represent the same element in $H^p(M,\rz)$. Let $\phi \in 
\mathcal{H}^p(Y)$ be the unique harmonic representative of the 
cohomology class of $\psi|_Y$. Then
 \begin{gather}
  \| F \|^2 \leq \| \psi \|^2_{L^2\Lambda^p(M)} + \frac{1}{\mu_1} \|\psi|_Y - 
\phi\|^2_{L^2\Lambda^p(Y)},
 \end{gather}
where $\mu_1^2$ is the smallest positive eigenvalue of $\Delta_Y$. In particular
 \begin{gather}
  \| F \|^2 \leq C(n,p) \mathrm{Vol}_*(M) \| [F|_M] \|_\infty^2,
 \end{gather}
 where we define the effective volume 
 $\mathrm{Vol}_*(M)$ by
 \begin{gather}
  \mathrm{Vol}_*(M) = \mathrm{Vol}(M) + \frac{1}{\mu_1} \mathrm{Vol}(Y).
 \end{gather}
\end{lem}
\begin{proof}
Let $F \in \mathcal{H}^p_{ext,abs}(X)$. Let $\psi \in
\Lambda^p(M)$ be a closed form such that 
\begin{gather}\label{closed2}
 F|_M - \psi = d h,
\end{gather}
for some $h \in \Lambda^{p-1}(M)$. Denote  by $\phi$
the unique harmonic representative of the class of $\psi|_Y$. If $\chi_Z$
is the characteristic function of $Z \subset X$, then the norm of $F$
is by definition the $L^2$-norm of $F - \phi \chi_Z$. Therefore,
\begin{gather}
 \| \hat \psi - \phi \chi_Z \| ^2 = \| F - \phi \chi_Z + (\hat \psi -F) \|^2 =\\=
 \| F - \phi \chi_Z\|^2 + \|(\hat \psi -F)\|^2 + 2 \langle F-\phi \chi_Z,(\hat \psi -F)\rangle. 
\end{gather}
Now observe that by (\ref{closed2}) and the definition of $\hat \psi$, 
the expansion of $\hat \psi -F$ on $Z$ contains no harmonic form. Hence by
(\ref{estimation8}) and (\ref{extension1}) the restriction of $\hat\psi -F$ 
to $Z$ is exponentially decreasing. This implies
\[
\langle F-\phi \chi_Z,\hat \psi -F\rangle= \langle F,\hat \psi -F\rangle.
\]
By (\ref{closed2}) and (\ref{derivative}) we have $F-\hat\psi=\widehat{dh}=
d\hat h$, where the latter means the distributional derivative.
Since $F$ is closed and coclosed, it follows that
\begin{gather}
\langle F-\phi \chi_Z,F-\hat \psi\rangle=\langle F,\widehat{dh}\rangle=
\langle F,d\hat h\rangle=\langle\delta F,\hat h\rangle =0.
\end{gather}
Therefore we get
\begin{gather}
 \| \hat \psi - \phi \chi_Z \| ^2 = \| F - \phi \chi_Z + (\hat \psi -F) \|^2 \nonumber =\\=
 \| F - \phi \chi_Z\|^2 + \|(\hat \psi -F)\|^2\ge\|F\|^2.
\end{gather}
On the other hand using an expansion of the form (\ref{expans6}) with a basis $(\phi_i)_{i \in \mathbb{N}}$
of co-exact eigenforms on the boundary we have
\begin{gather} \nonumber
 \| \hat \psi - \phi \chi_Z \| ^2 = \| \psi \|^2_{L^2\Lambda^p(M)} + 
 \int_0^\infty \sum_{i=1}^\infty |a_i|^2 e^{-2 \mu_{\phi_i} u} \left( \| d' \phi_i\|^2 
+ \mu_{\phi_i}^{2}
 \|\phi_i\|^2 \right) du = \\ = \nonumber
 \| \psi \|^2_{L^2\Lambda^p(M)} + \sum_{i=1}^\infty \mu_{\phi_i}^{-1} 
|a_i|^2 \|d' \phi_i\|^2
 = \| \psi \|^2_{L^2\Lambda^p(M)} +  \| \Delta_Y^{-\frac{1}{4}} (\psi|_Y - \phi)
 \|^2_{L^2\Lambda^p(Y)} \leq \\ \leq  \| \psi \|^2_{L^2\Lambda^p(M)} 
+ \frac{1}{\mu_1} \|\psi|_Y-\phi\|^2.
\end{gather}
\end{proof}

\begin{lem} \label{hilfreichlemma2}
Let $\phi \in \mathcal{H}^p_{ext,rel}(X)$. 
Let $\psi \in\Lambda^p(M,Y)$ be a representative of the class
$[\psi] \in H^p(M,Y,\rz)$ which corresponds to $\phi$ with respect to the 
isomorphism {\rm (\ref{caniso})}. Then
 \begin{gather}
  \| \phi \|^2 \leq \| \psi \|^2_{L^2\Lambda^p(M)},
 \end{gather}
 and in particular
 \begin{gather}
  \| \phi \|^2 \leq C(n,p) \mathrm{Vol}(M) \| [\psi] \|_\infty^2.
 \end{gather}
\end{lem}
\begin{proof} 
Let $\phi \in \mathcal{H}^p_{ext,rel}(X)$.
Recall the definition of $R_c$ by (\ref{derham1}). Choose $\chi$ such that the 
support of $\phi-d(\chi(u\phi_n+\theta))$ is contained in $M$. Then
$[\phi-d(\chi(u\phi_n+\theta))]$ is the image of $\phi$ w.r.t. the isomorphism
(\ref{caniso}). Since $\psi\in\Lambda^p(M,Y)$ represents this cohomology class,
there is $\omega\in\Lambda^{p-1}(M,Y)$ such that
\begin{gather}\label{equ6}
\phi-d(\chi(u\phi_n+\theta))=\psi+d\omega.
\end{gather}
Let $\tilde \psi$ (resp. $\tilde\omega$) be the differential form on $X$ 
which is equal to $\psi$ (resp. $\omega$) on $M$ and $0$ on $Z$.  Then 
 \begin{gather}\label{ungleichungsstifter}
\begin{split}
\| \psi \|^2&= \| \tilde\psi - \phi+ \chi_Z du \wedge\phi_n + (\phi - \chi_Z du
\wedge \phi_n) \|^2
=\| \tilde\psi - \phi + \chi_Z du \wedge \phi_n\|^2\\
&\quad + \| \phi - \chi_Z du\wedge \phi_n\|^2 
+ 2 \langle \tilde\psi - \phi + \chi_Z du \wedge \phi_n, \phi - \chi_Z du 
\wedge\phi_n\rangle. 
\end{split}
\end{gather}
By the definition (\ref{innprod1}) of the norm in $\cH^p_{ext,rel}(X)$, we have
\[
\parallel \phi\parallel =\parallel \phi- \chi_Z du \wedge \phi_n\parallel_{L^2}.
\]
By (\ref{equ6}) we have 
\begin{gather*} 
\begin{split}
\langle\tilde\psi - \phi + \chi_Z du \wedge \phi_n, \phi- \chi_Z du \wedge 
\phi_n \rangle&=\langle\tilde\psi - \phi, 
\phi- \chi_Z du \wedge \phi_n \rangle\\
&= -\langle d\tilde\omega +d(\chi(u\phi_n+\theta)), 
\phi - \chi_Z du \wedge\phi_n\rangle. 
\end{split} 
\end{gather*}
Since $\omega\in\Lambda^{p-1}(M,Y)$, it follows from Green's formula that
\[
\langle d\tilde\omega, \phi - \chi_Z du \wedge\phi_n\rangle=
\int_M d\omega\wedge *\phi=\int_Y i^*_Y(\omega)\wedge i^*_Y(*\phi)=0.
\]
Similarly, by Green's formula and (\ref{restrict}) we have
\[
\begin{split}
\langle d(\chi(u\phi_n+\theta)), \phi - &\chi_Z du \wedge\phi_n\rangle=
\int_M d(\chi(u\phi_n+\theta))\wedge *\phi\\
&+
\int_Z d(\chi(u\phi_n+\theta))\wedge *(\phi - \chi_Z du \wedge\phi_n)\\
&=\int_Y\theta\wedge *(du\wedge \phi_n+d\theta)-\int_Y \theta\wedge * d\theta=0.
\end{split}
\]
Thus
\[
\langle\tilde\psi - \phi + \chi_Z du \wedge \phi_n, \phi- \chi_Z du \wedge 
\phi_n \rangle=0
\]
and by (\ref{ungleichungsstifter}) we get
\[
\| \psi \|^2 =\| \phi \|^2+\| \tilde\psi - \phi + \chi_Z du \wedge \phi_n\|^2
\ge \| \phi \|^2.
\] 
\end{proof}

\section{Estimates on the scattering matrix and stable systoles} \label{systoles}

We recall some notions from geometric measure theory.
Suppose $B$ is a compact oriented Riemannian manifold and let $A$
be a closed submanifold. In our case $A$ will be either the boundary of $B$ or
the empty set.
For $z \in H_p(B,A,\zz)$ let the minimal volume be defined as the
infimum of the volumes of all its representatives, i.e.
\begin{gather}
 \mathrm{vol}(z)=\inf \{ \sum_i |\alpha_i| \mathrm{Vol}(c_i) \mid
 z= \sum_i \alpha_i [c_i],\;\alpha_i\in\Z\}.
\end{gather}
where the infimum is over all Lipschitz continuous simplices $c_i$.
The stable norm $\| z \|_{st}$ of an element $z \in H_p(B,A,\rz)$
is defined similarly by
\begin{gather}
 \| z \|_{st}:=\inf  \{ \sum_i |\alpha_i| \mathrm{Vol}(c_i) \mid 
 z = \sum_i \alpha_i [c_i],\;\alpha_i\in\R \}.
\end{gather}
This defines indeed a norm (see \cite[9.6, 9.9]{FeFl60}, and 
\cite[\S 3]{Fed74},  
also \cite[Ch. 4C]{Gro99}, and \cite[5.1.6]{Fed69}).
The stable norm of an element in $z \in H_p(B,A,\zz)$ is by definition the
stable norm of its image in $H_p(B,A,\rz)$. Clearly,
\begin{gather}
 \|z\|_{st} \leq \mathrm{vol}(z)
\end{gather}
and equality does not hold in general. However, as shown by Federer 
(see \cite[\S 5]{Fed74}),
\begin{gather}
 \|z\|_{st}= \lim_{k \to \infty} \frac{1}{k} \mathrm{vol}(k z).
\end{gather}
Moreover, for $z \in H_{n-1}(B,A,\zz)$ we have equality, i.e. $\|z\|_{st} =
\mathrm{vol}(z)$.

By a general result in geometric measure theory 
the stable norm and the comass are dual to each other, i.e.,  
\begin{gather}
 \| z \|_{st}= \sup \{ |\phi(z)| \mid \phi \in H^p(B,A,\rz), \| \phi \|_\infty \leq 1 \}
\end{gather}
(see \cite[4.10]{Fed74}, and also \cite[4.35]{Gro99},  and \cite{AuB06} for a 
sketch of the proof in the case without boundary).

\subsection{Estimate of the $L^2$-norm on $H^*(Y)$}

Now we will apply this result to our problem in the case where 
$B=Y$ and $A= \emptyset$. We equip $H^p(Y)$ with the norm induced from 
$\cH^p(Y)$ by the de Rham isomorphism.
Suppose that $\omega$ is a $p$-form on $Y$.
Then using (\ref{estconst}) one obtains
\begin{gather} \label{inform}
 \| \omega \|_2^2 \leq C(n-1,p)\mathrm{Vol}(Y) \| \omega \|_\infty^2 ,
\end{gather}
Now let $\phi$ be a harmonic $p$-form representing an element 
$[\phi] \in H^p(Y,\rz)$. Using (\ref{inform}), we get 
\begin{gather*}
\begin{split}
 \| \phi \|_2 &= \sup \left\{ \int\phi \wedge \omega | \:\omega \in 
\Lambda^{n-p-1}(Y), d\omega=0, 
 \|\omega\|_2 \leq 1\right\} \\ 
&\geq C(n-1,p)^{-1/2} \mathrm{Vol}(Y)^{-1/2}
 \sup \left\{ \int \phi \wedge \omega | \:\omega \in \Lambda^{n-p-1}(Y), 
d\omega=0, \|\omega\|_\infty \leq 1\right\}\\ 
&=C(n-1,p)^{-1/2} \mathrm{Vol}(Y)^{-1/2} 
 \sup \left\{ \langle [\phi] \cup \alpha, [Y] \rangle \mid \alpha \in 
H^{n-p-1}(Y,\rz), \| \alpha \|_\infty \leq 1 \right\}
\end{split}
\end{gather*}
Since the stable norm is dual to the comass norm, we finally get
\begin{gather}
 \| \phi \|_2 \geq C(n-1,p)^{-1/2}  \mathrm{Vol}(Y)^{-1/2} \| [Y] \cap 
\phi \|_{st},
\end{gather}
for any $\phi \in \cH^p(Y,\rz)$. 

On the other hand the inequality (\ref{inform}) may also be used 
directly. Since $\| \phi \|_2$ is the infimum of the $L^2$-norms of all  
representatives of the cohomology class $[\phi]$, we have
\begin{pro} \label{randabsch}
The Hilbert space norm on $H^p(Y,\rz)$, induced from the harmonic forms, 
satisfies
\begin{gather}
 C(n-1,p)^{-1/2} \mathrm{Vol}(Y)^{-1/2} \| [Y] \cap \phi \|_{st} \leq 
\|\phi\|_2 \leq
 C(n-1,p)^{1/2}  \mathrm{Vol}(Y)^{1/2} \|\phi \|_\infty.
\end{gather}
\end{pro}

\subsection{Estimate of the norms on $\mathcal{H}_{ext}(X)$}

We can use Lemmas \ref{hilfreichlemma} and \ref{hilfreichlemma2} 
in the same way as above to get similar estimates of the norms of
elements in $\mathcal{H}^p_{ext,rel}(X)$. First note that $*$ induces an
isomorphism
\[
*\colon \cH_{ext,rel}^p(X)\to \cH_{ext,abs}^{n-p}(X).
\]
Let $(\cdot,\cdot)_{ext}$ be the pairing (\ref{pairing}). Then we have
\[
\langle \phi,\psi\rangle =(\phi,*\psi)_{ext},\quad \phi,\psi\in 
\cH^p_{ext,rel}(X).
\]
Let $F \in\mathcal{H}^p_{ext,rel}(X)$. Then we get
\begin{gather}\label{norm10}
\begin{split}
 \| F \| &= \sup \{ \langle F, \omega \rangle \mid \omega \in 
\mathcal{H}^{p}_{ext,rel}(X), \|\omega\| \leq 1\}\\ 
&=\sup \left\{(F,\psi)_{ext} \mid \psi\in \mathcal{H}^{n-p}_{ext,abs}(X),\;  
\|\psi\|\le 1\right\}.
\end{split}
\end{gather}
Choose a representative $\phi$ of
the cohomology class $R_c(F)\in H_c^p(X)$ with $\supp(\phi)\subset M$. Then
$\eta_{rel}(F)=[\phi]$. Since the diagramm (\ref{diagram}) commutes, we get
\[
\langle \phi,\psi\rangle=(R_c(\phi),R(*\psi))=([\phi],[*\psi|_M]), \quad 
\phi,\psi\in H^p_{ext,rel}(X).
\]
Using Lemma \ref{hilfreichlemma}, it follows from (\ref{norm10}) that
\[
\begin{split}
\| F \|_2&=\sup\{ ([\phi],[\psi|_M]) \mid \psi\in 
\mathcal{H}^{n-p}_{ext,abs}(X),\,  \|\psi\|\le 1\}\\
&\ge C(n,p)^{-1/2} \mathrm{Vol}_*(M)^{-1/2} 
 \sup \{ \langle [\phi] \cup \alpha, [M] \rangle \mid \alpha \in H^{n-p}(M,\rz), \| \alpha \|_\infty \leq 1 \}.
\end{split}
\]
Since the stable norm is dual to the comass norm, we finally get
\begin{gather}
 \| F \| \geq C(n,p)^{-1/2}  \mathrm{Vol}_*(M)^{-1/2} \| [M] \cap 
[\phi] \|_{st},
\end{gather}

Note that the Poincare-Lefschetz dual $[M] \cap [\phi]$
of $[\phi]$ is in $H_p(M,\rz)$ and its stable norm equals its stable norm as 
an 
element of $H_p(M,\rz)$. Lemma \ref{hilfreichlemma2} gives an upper bound
for $\|F\|$ and we get 

\begin{pro}
 Let $F\in\mathcal{H}^p_{ext,rel}(X)$ and $\phi=\eta_{rel}(F)\in H^p(M,Y,\rz)$. 
Then we have
\begin{gather}
 C(n,p)^{-1/2} \mathrm{Vol}_*(M)^{-1/2} \| [M] \cap \phi \|_{st} \leq 
\|F\| \leq
 C(n,p)^{1/2}  \mathrm{Vol}(M)^{1/2} \|\phi \|_\infty.
 \end{gather}
\end{pro}
Using Theorem \ref{main3}, we obtain Theorem \ref{main00}.

If we denote by $\iota$ the inclusion map $Y \to M$ and by $\iota_*: H_p(Y,\rz) \to
H_p(M,\rz)$ the induced map in homology, then
\begin{gather} \label{anders}
 [M] \cap \partial \phi = \iota_*\left([Y] \cap \phi \right).
\end{gather} 
That is $[M] \cap \partial \phi$ coincides with the image of the
Poincare dual of the class $\phi$ in $H_{n-p-1}(M,\rz)$. The
stable norm of $\partial \phi$ can be calculated in many cases explicitly
in terms of geometric data using the fact that the stable mass is dual to the
comass norm. In order to make statements about the spectrum of the map $T(0)$
one can combine these estimates with the estimates of the
$L^2$-norms of the cohomology classes on $Y$ (see Proposition \ref{randabsch}).

\section{Examples} \label{examples}

\subsection{The scattering length for functions}
Note that the extended harmonic functions are exactly the constant
functions. Therefore,  $\cH^0(Y) \cap \ker \partial$ is spanned by the function equal
to $1$ on $Y$. Thus, the $+1$ eigenspace to $S_0(0)$ is spanned by $1$.
Moreover, since $S$ anticommutes with the Hodge $*$
operator we have $S_{n-1}(0) *1= -*1$. 
By Stokes formula
\begin{gather}
 (\partial [*1])[M]=[* 1][Y]
\end{gather}
and therefore
\begin{gather}
 \partial [*1] = \frac{\mathrm{Vol}(Y)}{\mathrm{Vol}(M)} [*_M 1]
\end{gather}
and by the above equation we immediately obtain
\begin{gather}
 T_{n-1}(0) |_{(\ker \partial)^\perp}= 2 \frac{\mathrm{Vol}(M)}{\mathrm{Vol}(Y)}.
\end{gather}
This in turn implies that
\begin{gather}
 T_0(0)|_{\ker \partial} =2  \frac{\mathrm{Vol}(M)}{\mathrm{Vol}(Y)}.
\end{gather}

\subsection{$Y$ has only one connected component}
In this case $\mathcal{H}^0(Y)$
consists of the constant functions only. By the above we have
\begin{gather}
 T_{n-1}(0)= 2 \frac{\mathrm{Vol}(M)}{\mathrm{Vol}(Y)}.
\end{gather}
This in turn implies that
\begin{gather}
 T_0(0)= 2 \frac{\mathrm{Vol}(M)}{\mathrm{Vol}(Y)}.
\end{gather}
Note that if $Y=S^{n-1}$  the only non-vanishing cohomology
groups are $H^0(Y)$ and $H^{n-1}(Y)$. Thus, in this case the above formulas determine
$T(0)$ completely.

\subsection{$Y=Y_1 \cup Y_2$ has two boundary components} \label{scatterex1}

Now we have $\cH^0(Y) \cong \rz^2$ and therefore on $\cH^0(Y)$ is a direct sum
of the one dimensional spaces $\ker{\partial}$ and $\ker{\partial}^{\perp}$.
Under the splitting $\cH^0(Y) = \ker{\partial} \oplus \ker{\partial}^{\perp}$
the operator $T_0(0)$ is of the form
\begin{gather}
 T_0(0) = \left ( \begin{matrix} t_1 & 0 \\ 0 & t_2 \end{matrix} \right), 
\end{gather}
where we have already seen, that 
\begin{gather}
 t_1 = 2 \frac{\mathrm{Vol}(M)}{\mathrm{Vol}(Y)}.
\end{gather}
Our formula now allows us to give an estimate of $t_2$ in purely geometric 
terms.
The harmonic functions in $\ker{\partial}$
are multiples of the function $1$ on $Y$. The complement 
$\ker{\partial}^\perp$ has dimension $1$ and is spanned by the function
\begin{gather}
 \chi(x) = \left \{ \begin{matrix}  \mathrm{Vol}(Y_2) \: \mathrm{for}\; x 
\in Y_1 \\
 -\mathrm{Vol}(Y_1) \: \mathrm{for}\; x \in Y_2  \end{matrix}  \right.
\end{gather}
Clearly,
\begin{gather}
 \| \chi \|_2^2 = \mathrm{Vol}(Y_2)^2  \mathrm{Vol}(Y_1) + 
\mathrm{Vol}(Y_1)^2  \mathrm{Vol}(Y_2). 
\end{gather}

It remains to estimate the $L^2$-norm of the class $\partial \chi$ in $
H^1(M,Y,\rz)$.
The comass norm of $\partial \chi$ may be calculated using the duality between
the stable norm and the comass norm. 
Let $c$ be a Lipschitz continuous chain
whose boundary is homologuous to $y_2 - y_1$, where
$y_1 \in Y_1$ and $y_2 \in Y_2$.
Then $c$  defines a relative cycle $[c]$ in $H_1(M,Y,\rz)$ and we have
\begin{gather}
 \partial \chi ([c]) = \int_c \partial \chi 
= \mathrm{Vol}(Y_1) +\mathrm{Vol}(Y_2)=\mathrm{Vol}(Y). 
\end{gather}
We observe that the cycle $c$ can be written as a linear combination of curves
which are either closed or whose end points are in the boundary. This implies
\begin{gather}
 \| [c] \|_{st} \geq \mathrm{dist}(Y_1,Y_2),
\end{gather}
with equality if for example $[c]$ is in the same homology class as a shortest curve 
connecting $Y_1$ with $Y_2$.

Duality implies
\begin{gather}
 \| \partial \chi \|_\infty = \mathrm{dist}(Y_1,Y_2)^{-1} \mathrm{Vol}(Y), 
\end{gather}
and therefore
\begin{gather}
 \| \partial \chi \|_2 \leq  \mathrm{Vol}(M)^{1/2} \mathrm{dist}(Y_1,Y_2)^{-1} 
\mathrm{Vol}(Y).
\end{gather}

To get the estimate from above we need to look at 
$\| [M] \cap \partial \chi \|_{st}$, 
that is the stable norm of the Poincare dual
of the class $\partial \chi$. By Equation (\ref{anders}) we have $[M] \cap
\partial \chi = \iota_*([Y] \cap \chi)$. Of course,
\begin{gather}
 \iota_*([Y] \cap \chi) = \iota_*(\mathrm{Vol}(Y_2)[Y_1] - \mathrm{Vol}(Y_1)[Y_2]) =
 \mathrm{Vol}(Y) \iota_*([Y_1]),
\end{gather}
where we have used that $\iota_*([Y_1])+\iota_*([Y_2])=0$. Therefore,
\begin{gather}
 \| M \cap \partial \chi \|_{st} = \mathrm{Vol}(Y) \|\iota_*([Y_1])\|_{st}.
\end{gather}

Combining these two estimates we obtain
\begin{gather}
 C_2 \leq t_2 \leq C_1
\end{gather}
with
\begin{gather}
 C_1= 2 \mathrm{Vol}_*(M) \frac{\mathrm{Vol}(Y_1) \mathrm{Vol}(Y_2)}
 {\|\iota_*([Y_1])\|_{st}^2 (\mathrm{Vol}(Y_1)+ \mathrm{Vol}(Y_2))},\\
 C_2= 2 \mathrm{Vol}(M)^{-1} \frac{\mathrm{dist}(Y_1,Y_2)^2 \mathrm{Vol}(Y_1)\mathrm{Vol}(Y_2)}
 {\mathrm{Vol}(Y_1)+ \mathrm{Vol}(Y_2)}.
\end{gather}

Note that with respect to this basis, the scattering matrix at zero has 
the form
\begin{gather}
 S_0(0) = \left ( \begin{matrix} 1 & 0 \\ 0 & -1 \end{matrix} \right). 
\end{gather}

In order to interpret this in terms of reflection and transmission coefficients
we choose another basis $(\chi_1,\chi_2)$ of $\cH^0(Y)$, where $\chi_i$
is constant equal to $1$ on $Y_i$ and equal to zero on the other boundary 
component. In this basis we get
\begin{gather}
 S_0(0) = \frac{1}{\mathrm{Vol}(Y_1)+\mathrm{Vol}(Y_2)}\left ( \begin{matrix} \mathrm{Vol}(Y_1)-\mathrm{Vol}(Y_2) & 
 2 \mathrm{Vol}(Y_2) \\ 2 \mathrm{Vol}(Y_1) 
 &  \mathrm{Vol}(Y_2)-\mathrm{Vol}(Y_1)\end{matrix} \right). 
\end{gather}
 The reflection coefficient $r_{11}$ and the transmission coefficient $r_{12}$
 of a wave of low energy that comes in at the end $Y_1$ and is scattered in $M$
 therefore are $r_{11}=\frac{\mathrm{Vol}(Y_1)-\mathrm{Vol}(Y_2)}{\mathrm{Vol}(Y_1)+\mathrm{Vol}(Y_2)}$ and
 $r_{12}=\frac{2 \mathrm{Vol}(Y_1)}{\mathrm{Vol}(Y_1)+\mathrm{Vol}(Y_2)}$.
 Transforming $T_0(0)$ into this basis gives
 \begin{gather}
 T_0(0) = \frac{1}{\mathrm{Vol}(Y_1)+\mathrm{Vol}(Y_2)}\left ( \begin{matrix} t_1 \mathrm{Vol}(Y_1) + t_2 \mathrm{Vol}(Y_2) & 
 (t_1-t_2) \mathrm{Vol}(Y_2) \\ (t_1-t_2) \mathrm{Vol}(Y_1) 
 &  t_1 \mathrm{Vol}(Y_2)+ t_2 \mathrm{Vol}(Y_1)\end{matrix} \right). 
\end{gather}
As a remark we would like to add here that in physics the scattering matrix for
one dimensional scattering problems has the transmission coefficients on the diagonal
and the reflection coefficients off the diagonal. Our notation differs here as we consider
the operator with Neumann boundary conditions at each end as the unperturbed operator.

\subsection{The full-torus} \label{scatterex2}

\begin{figure}
  \centerline{\includegraphics*[width=7cm]{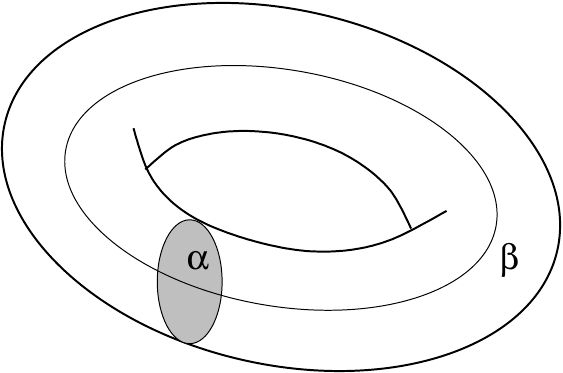}}
  \caption{The cycles $\alpha \in H_2(M,Y)$ and $\beta \in H_1(M)$ on $M$.} \label{fiig02}
\end{figure}

Let $M$ be the full torus $D \times S^{1}$ with boundary $Y=T^2=S^1 \times S^1$.
We view both $D$ and $S^1$ as subsets of the complex planes and use
coordinates $z=r e^{\I y}$ on $D$ and $z=e^{\I x}$ on $S^1$. We assume that
we are given a metric on $M$ which has product structure in a small
neighborhood of the boundary and that the metric on the boundary is equal to
the product metric
\begin{gather}
 \ell_1^2 dx^2 + \ell_2^2 dy^2,
\end{gather} 
with positive real numbers $\ell_1,\ell_2$. Then, the volume is 
$\mathrm{Vol}(Y)=4 \pi^2 \ell_1 \ell_2$. Moreover, 
$H^1(M,\rz) \cong \rz$ and it is gerated by the class of the
one form $dx$. The group $H^1(Y,\rz)$ is isomorphic to $\rz^2$ and is
generated by the classes of the two harmonic $1$-forms $dx$ and $dy$.
The $L^2$-norms of these forms is easily calculated.
\begin{gather}
 || dx ||^2_{L^2} = \int_Y dx \wedge * dx = 4 \pi^2 \frac{\ell_2}{\ell_1},\\ 
 || dy ||^2_{L^2} = \int_Y dy \wedge * dy = 4 \pi^2 \frac{\ell_1}{\ell_2}.
\end{gather}
and, moreover, $dx$ and $dy$ are orthogonal to each other.
The restriction of the form $dx$ on $M$ to $Y$ is the form $dx$ regarded as a
form on $Y$. Therefore, the kernel  of the connecting homomorphism $\partial$ is spanned by $[dx]$.
Hence, $(\ker\partial)^{\perp}$ is spanned by  $[dy]$.
Since the Hodge star operator commutes with $T_2(0)$ the map $T_2(0)$
has the form
\begin{gather}
 T_2(0) = \left ( \begin{matrix} t & 0 \\ 0 & t \end{matrix} \right), 
\end{gather}
 with respect to the decomposition $H^1(T^2,\rz)=\ker\partial \oplus
 (\ker\partial)^{\perp}$.
In order to get an estimate for $t$ we need to calculate $\| \partial [dy]\|_\infty$
and $\|[M] \cap \partial [dy]|_{st}$.
The stable norm of $\| \partial [dy]\|_\infty$ can again be calculated by
duality. For this note that $H_2(M,Y,\rz)\cong \rz$ by Alexander duality
and it is generated by the relative cycle $\alpha=D \times \{1\} \subset M$.
Now, of course
\begin{gather}
 \langle \partial [dy],[\alpha] \rangle= \int_{\partial D \times \{1\}} dy = 2 \pi.
\end{gather}
Duality of comass norm and stable norm now implies that
\begin{gather}
 \|\partial [dy]\|_{\infty}=2 \pi \|[D \times \{1\}]\|_{st}^{-1}.
\end{gather}
The class $[M] \cap \partial [dy]$ is given by $-2 \pi \iota_*(\{1\} \times S^{1})$.
So $\frac{1}{2 \pi}\|[M] \cap \partial [dy]\|_{st}$ is equal to the infimum of the lengths
of all representatives of the cycle $\beta=\{1\} \times S^{1}$ in $H_1(M,\rz)$.
The geometric picture is described in Fig. \ref{fiig02}.
Theorem \ref{main00} now gives an estimate of $t$:
\begin{gather}
 2 \frac{\ell_1}{\ell_2} \frac{\| \alpha \|_{st}^2}{\mathrm{Vol}(M)}
 \leq t \leq 2 \frac{\ell_1}{\ell_2} \frac{\mathrm{Vol}_*(M) }{\|\beta\|_{st}^2}.
\end{gather}
Since $\mu_1=\min \{\ell_1^{-1},\ell_2^{-1}\}$ we have
$\mathrm{Vol}_*(M)=\mathrm{Vol}(M) + 4 \pi^2 \ell_1 \ell_2 \max \{\ell_1,\ell_2\}$.

\appendix

\section{Dynamical approach and spectral decomposition} \label{app2}

In this appendix we discuss the relation between the stationary and the 
dynamical approach to scattering theory and we establish the Eisenbud-Wigner
formula for manifolds with cylindrical ends. For details concerning scattering
theory we refer to \cite{Ya92}. For original papers on the Eisenbud-Wigner formula
see \cite{Eis48} and \cite{Wi55}. In order to simplify the relation
of the scattering length to the time-delay operator we consider scattering
theory for the square root of the Laplace Beltrami operator, 
i.e. we consider scattering theory and the time-delay in relativistic quantum mechanics.

Let $\overline{\Delta_p}$ be the closure in $L^2$ of the operator $\Delta_p$ with domain
$\Lambda^p_c(X)$. Denote by $\Delta_{p,0}$ the self-adjoint
extension of $\Delta_p$ with respect to Neumann boundary conditions along $Y$. 
Note that
$\Delta_{p,0}$ is the self-adjoint operator associated to the quadratic form
\begin{gather*}
 \phi \mapsto \int_M q_x(\phi) dx + \int_Z q_x(\phi),\\
 q_x(\phi) = d \phi(x) \wedge \overline{* d \phi(x)} + \delta \phi(x) \wedge 
\overline{* \delta \phi(x)} 
\end{gather*}
on $H^1(M,\Lambda^p T^* M) \oplus H^1(Z,\Lambda^p T^* Z) \subset 
L^2\Lambda^p(X)$. With respect to the decomposition
\[
L^2\Lambda^p(X)=L^2\Lambda^p(M)\oplus L^2\Lambda^p(Z)
\]
we have
\[
\Delta_{p,0}=\Delta_{p,M}\oplus\Delta_{p,Z} ,
\]
where $\Delta_{p,M}$ and $\Delta_{p,Z}$ are the corresponding self-adjoint extensions of
$\Delta_p|_M$ and $\Delta_p|_Z$, respectively, with Neumann boundary conditions
imposed. Let $H$ be the square root of $\overline{\Delta_p}$ defined by spectral calculus
and similarly define $H_0=(\Delta_{p,0})^{1/2}$, $H_M=\Delta_{p,M}^{1/2}$,
and $H_Z=\Delta_{p,Z}^{1/2}$.
Since $M$ is compact, $H_M$ has purely discrete spectrum. The spectral
resolution of $H_Z$ can be determined, using separation of variables. The
spectrum is absolutely continuous. A complete set of generalized eigensections 
is given by
\begin{gather}
  F_0(\phi,\lambda)=\left( e^{+\I \sqrt{\lambda^2-\mu_\phi^2} u} + e^{-\I
     \sqrt{\lambda^2-\mu_\phi^2} u}\right) j_p(\phi),
\end{gather}
where $\phi$ runs through an orthonormal basis of  eigensection of the 
operator $\Delta_p' \oplus\Delta_{p-1}'$ with eigenvalue $\mu_\phi^2$.

Denote by $P_{ac}$ (resp. $P_{ac}^0$) the orthogonal projection onto
the absolutely continuous subspace of $H$ (resp. $H_0$). 
Denote by  $H_{ac}$ and $H_{0,ac}$ the restrictions of $H$ and $H_0$, 
respectively, to the absolutely continuous subspaces. Furthermore, for an open 
interval
$(a,b)$ we denote by $P_{ac}(a,b)$ (resp. $P_{ac}^0(a,b)$) the projection onto 
the continuous part of the spectral subspace of the interval. As explained 
above, the absolutely continuous subspace for $H_0$ is $L^2\Lambda^p(Z)$. Thus
$P_{ac}^{0}$ is the orthogonal projection of $L^2\Lambda^p(X)$ onto the 
subspace $L^2\Lambda^p(Z)$. 

By \cite[Th\'eor\`em 3.6]{Gui89} and the Birman-Kato invariance principle
the wave operators
\begin{gather}
 W_\pm = \mathrm{s}-\lim_{t \to \pm \infty} e^{\I t H} e^{-\I t H_0} P_{ac}^0
\end{gather}
exist and are complete. This means that the strong limit exists and and the
operators $W_\pm$ define isometries of $P_{ac}^0L^2\Lambda^p(X)$ onto 
$P_{ac}L^2\Lambda^p(X)$, intertwining $H_{0,ac}$ and $H_{ac}$. 
In this context the scattering operator $S$ is defined as
\begin{gather}
 \mathcal{S} = W_+^* W_-.
\end{gather}
This is a unitary operator in $P_{ac}^0 L^2\Lambda^p(X)$ which commutes with 
$H_{0,ac}$. Let $\sigma_0=\sigma_{ac}(H_0)$ be the absolutely continuous 
spectrum of $H_0$. It equals $[\mu,\infty)$ with $\mu\ge0$. Let 
$\left\{E_0(\lambda)\right\}_{\lambda\in\sigma_0}$ be
the spectral family of $H_{0,ac}$. Since $S$ commutes with $H_{0,ac}$, we have
\[
S=\int_{\sigma_0} S(\lambda)\;dE_0(\lambda),
\]
where $S(\lambda)=S(\lambda;H,H_0)$ acts in the finite-dimensional Hilbert 
space 
\begin{gather} 
 \mathcal{H}(\lambda) = \bigoplus_{\mu \leq \lambda^2} 
\mathrm{Eig}_\mu(\Delta'_p) \oplus \mathrm{Eig}_\mu(\Delta'_{p-1}).
\end{gather}

For the purposes of this paper we need only to investigate the structure of 
the continuous spectrum of $H$ in the interval $[0,\mu_1)$ in the situations
when $\mu=0$. 
The generalized eigensections $F(\phi,\lambda)$
constructed in Theorem \ref{main1} are well defined as distributions
on $(0,\mu_1)$ with values in $L^2\Lambda^p(X)$. This follows easily
from the properties of the Fourier transform and the expansion of
$F(\phi,\lambda)$ on $Z$. Of course, in the weak topology on the space of
distributions we have
\begin{gather}
 \langle F(\phi,\lambda), F(\psi,\lambda') \rangle = 
 \lim_{a \to \infty} \int_{M_a} F(\phi,\lambda)(x) \wedge \overline{* F(\psi,\lambda')(x)}.
\end{gather}
The smooth function
\begin{gather}
 \int_{M_a} F(\phi,\lambda)(x) \wedge \overline{* F(\psi,\lambda')(x)}
\end{gather}
can be determined as the continuous extension of
\begin{gather}
 \frac{1}{\lambda^2-(\lambda')^2}\int_{M_a} 
\Delta_p F(\phi,\lambda) \wedge \overline{*
 F(\psi,\lambda')(x)} - F(\phi,\lambda) \wedge 
\overline{* \Delta_p F(\psi,\lambda')(x)}.
\end{gather}
This integral can be simplified using Green's formula and the limit
$a \to \infty$ can be explicitly evaluated (see Equ. (\ref{parequa})). A simple exercise in distribution
theory shows that indeed we have the orthogonality relations
\begin{gather} 
\langle F(\phi,\lambda), F(\psi,\lambda') \rangle = 2 \pi \langle \phi,\psi \rangle
 \delta(\lambda-\lambda') 
\end{gather}
 as distributions on $(0,\mu_1) \times (0,\mu_1)$. Moreover, in the same
 way as in \cite{Gui89}, Theorem 6.2, one shows that
\begin{gather}
 W_+ F_0(\phi,\lambda) = F(C^*(\lambda)\phi,\lambda),\\
 W_- F_0(\phi,\lambda) = F(\phi,\lambda).
\end{gather}
Therefore, in the distributional sense for $\lambda \in (0,\mu_1)$
\begin{gather}
  \mathcal{S} F_0(\phi,\lambda) = F_0(C_p(\lambda),\lambda) 
\end{gather}
which shows that the dynamical and the stationary scattering matrices coincide
\begin{gather}
 \mathcal{S}(\lambda)=C_p(\lambda).
\end{gather}

The time delay operator $\T$ is defined in the following way. If 
$\phi \in P_{ac} L^2\Lambda^p(X)$
then, according to the laws of quantum mechanics, the probability of finding 
the particle
with wave function $\phi$ in $M_a$ at time $t$ is given by
\begin{gather}
 \int_{M_a} \| e^{-\I Ht} \phi \|_x^2 dx = \| \chi_{M_a} e^{-\I H t} \phi \|^2.
\end{gather}
The total time spent in $M_a$ is then given by
\begin{gather}
 \int_{-\infty}^{\infty}\|\chi_{M_a} e^{-\I H t} \phi \|^2 dt.
\end{gather}
This expression is not necessarily finite for all $\phi$. Now, according to 
scattering theory,
for $\phi \in P_{ac}^0 L^2\Lambda^p(X)$ the states $e^{-\I t H}W_- \phi$ and 
$e^{-\I H_0 t}\phi$ 
are asymptotically the same for $t \to - \infty$. Thus, the time excess due to 
the interaction (the presence of M) is
\begin{gather}
\int_{-\infty}^{\infty} \left( \|\chi_{M_a} e^{-\I H t} W_- \phi \|^2
- \|\chi_{M_a} e^{-\I H_0 t}  \phi \|^2 \right)dt.
 \end{gather}
The Eisenbud-Wigner time-delay operator $\T$ is formally defined by
\begin{gather}
 \langle \phi, \T \phi \rangle = 
 \lim_{a \to \infty} \int_{-\infty}^{\infty} \left( \|\chi_{M_a} e^{-\I H t} W_- 
 P_{ac}^0 \phi \|^2- \|\chi_{M_a} e^{-\I H_0 t} P_{ac}^0  \phi \|^2 \right)dt.
\end{gather}
 In many situations in potential scattering it can be shown that the above 
defines
a closable quadratic form and $\T$ is a self-adjoint operator that commutes 
with $H_0$ and the Eisenbud-Wigner formula
 \begin{gather}
  \T = \int_{\sigma_{ac}(H_0)} \T(\lambda)\; d E_0(\lambda),\\
  \T(\lambda) = - \I \mathcal{S}(\lambda)^{-1} 
\frac{d \mathcal{S}}{d\lambda}(\lambda)
 \end{gather}
 holds. Since we are only interested  in the spectrum near $0$, we prove this 
formula  only for elements in $P^0_{ac}(0,\mu_1)$.
\begin{pro}
Suppose that $g \in C^{\infty}_0(0,\mu_1)$, $\phi \in \mathrm{ker} 
\Delta_p' \oplus \Delta_{p-1}'$. 
 Let $F_0(\phi,g):=\int_\R F_0(\phi,\lambda) g(\lambda) d \lambda$. Then,
 \begin{gather*}
 \int_{-\infty}^{\infty}\|\chi_{M_a} e^{-\I H t} W_- F_0(\phi,g) \|^2 dt 
 < \infty,\\
 \int_{-\infty}^{\infty}\|\chi_{M_a} e^{-\I H_0 t}  F_0(\phi,g) \|^2 dt < \infty.
 \end{gather*}
 Moreover, 
  \begin{gather} \label{bloedequ}
    \langle F_0(\phi,g),  \T F_0(\phi,g) \rangle 
= 2 \pi \int_{\sigma_0} \T(\lambda) g(\lambda)^2 \;d E_0(\lambda),
  \end{gather}
  where 
  \begin{gather}
   \T(\lambda)=- \I \mathcal{S}(\lambda)^{-1} \mathcal{S}'(\lambda).
  \end{gather}
  Hence, $\T$ is a self-adjoint operator on $P^0_{ac}(0,\mu_1)L^2\Lambda^p(X)$
  which has the form
  \begin{gather}
  \T = \int_{\sigma_0} \T(\lambda)\; d E_0(\lambda)
  \end{gather}
 with respect to the spectral family 
$\left\{E_0(\lambda)\right\}_{\lambda\in\sigma_0}$ of $H_{0,ac}$.
 \end{pro}
 \begin{proof}
   We can do this calculation in the explicit spectral decomposition.
  \begin{gather} \nonumber
   \int_{-\infty}^{\infty}\|\chi_{M_a} e^{-\I H t}  F_0(\phi,g) \|^2 d \lambda = \\=
   \int_{-\infty}^{\infty} \int \int \int_{M_a} e^{-\I(\lambda-\lambda') t} \langle F_0(\phi,\lambda), F_0(\phi,\lambda') \rangle_x
   g(\lambda) g(\lambda') dx d \lambda d \lambda' dt = \\ =
   2 \pi \int  \langle F_0(\phi,\lambda), \chi_{M_a} F_0(\phi,\lambda) \rangle g(\lambda)^2 d \lambda. \nonumber
  \end{gather}
  and similarly,
  \begin{gather} \nonumber
   \int_{-\infty}^{\infty}\|\chi_{M_a} e^{-\I H t}  F(\phi,g) \|^2 d \lambda = \\=
   2 \pi \int \langle \langle F(\phi,\lambda), \chi_{M_a} F(\phi,\lambda) \rangle g(\lambda)^2 d \lambda.
  \end{gather}
   In the limit $a \to \infty$ both integrals can be computed up to exponentially small
   errors since we have the expansions on $Z$:
   \begin{gather}
     F(\phi,\lambda)|_Z = e^{-\I \lambda u} j_p(\phi) + e^{+\I \lambda u} j_p(S(\lambda) \phi) + R,\\
     F_0(\phi,\lambda)|_Z = e^{-\I \lambda u} j_p(\phi) + e^{+\I \lambda u} j_p(\phi).
   \end{gather}
   Therefore, integration by parts yields
   \begin{gather} \nonumber
    \langle F_0(\phi,\lambda), \chi_{M_a} F_0(\phi',\lambda') \rangle = \\=
    \frac{1}{\lambda^2-(\lambda')^2} \int_{M_a} \langle \Delta_p F_0(\phi,\lambda), F_0(\phi',\lambda') \rangle_x -
    \langle F_0(\phi,\lambda), \Delta_p F_0(\phi',\lambda') \rangle_x = \\ =
    \frac{2}{\lambda-\lambda'} \sin\left((\lambda-\lambda') a\right) \langle \phi,\phi' \rangle +
    \frac{2}{\lambda+\lambda'} \sin\left((\lambda+\lambda') a\right) \langle \phi,\phi' \rangle. \nonumber
   \end{gather}
   Similarly, we have
   \begin{gather}  \nonumber
    \langle F(\phi,\lambda), \chi_{M_a} F(\phi',\lambda') \rangle = \\= \nonumber
    \frac{1}{\lambda^2-(\lambda')^2} \int_{M_a} \langle \Delta_p F(\phi,\lambda), F(\phi',\lambda') \rangle_x -
    \langle F(\phi,\lambda), \Delta_p F(\phi',\lambda') \rangle_x = \\ = \label{parequa}
    \frac{2}{\lambda-\lambda'} \sin\left((\lambda-\lambda') a\right) \langle \phi,\phi' \rangle +	
    \left(\frac{\I}{\lambda-\lambda'} e^{\I(\lambda-\lambda')a} \right) 
    \langle \phi, (\id - \mathcal{S}^*(\lambda) \mathcal{S}(\lambda') )\phi' \rangle -\\-
    \frac{\I}{\lambda + \lambda'}\left( e^{\I (\lambda+\lambda')a} \langle \phi, 
    \mathcal{S}^*(\lambda) \phi'\rangle -
     e^{-\I (\lambda+\lambda')a} \langle \phi, \mathcal{S}(\lambda') \phi'\rangle  \right) 
    + O(e^{-\mu_1 a}). \nonumber
   \end{gather}
   Therefore,
   \begin{gather*}
    \langle F(\phi,\lambda), \chi_{M_a} F(\phi',\lambda') \rangle - 
    \langle F_0(\phi,\lambda), \chi_{M_a} F_0(\phi',\lambda') \rangle= \\=
    \left(\frac{\I}{\lambda-\lambda'} e^{\I(\lambda-\lambda')a}\right) \langle \phi, 
    (\id - \mathcal{S}^*(\lambda) \mathcal{S}(\lambda')) \phi' \rangle -
     \frac{2}{\lambda+\lambda'} \sin\left((\lambda+\lambda') a\right) \langle \phi,\phi' \rangle-\\-
     \frac{\I}{\lambda + \lambda'}\left( e^{\I (\lambda+\lambda')a} \langle \phi, \mathcal{S}^*(\lambda) \phi'\rangle -
     e^{-\I (\lambda+\lambda')a} \langle \phi, \mathcal{S}(\lambda') \phi'\rangle  \right)+
    O(e^{-\mu_1 a}),
   \end{gather*}
   and one obtains
   \begin{gather}
     \lim_{a \to \infty} \lim_{\lambda' \to \lambda}\langle F(\phi,\lambda), \chi_{M_a} F(\phi',\lambda) \rangle - 
    \langle F_0(\phi,\lambda), \chi_{M_a} F_0(\phi',\lambda) \rangle =
     \langle \phi,\T(\lambda) \phi' \rangle,
   \end{gather}
   where the second limit is in the distributional sense and 
   \begin{gather}
    \T(\lambda)=-\I \mathcal{S}(\lambda)^{-1} \mathcal{S}'(\lambda).
   \end{gather}
 \end{proof}

By the properties of the scattering matrix  $\mathcal{T}(\lambda)$ commutes
with the Hodge star operator and leaves the summands in $\mathrm{ker} \Delta_p' \oplus \Delta_{p-1}'$
invariant (see (\ref{invsummands})). We therefore have
\begin{gather} 
  \T_p(\lambda)=\left( \begin{matrix} T_p(\lambda) && 0\\ 0 &&  T_{p-1}(\lambda) \end{matrix} \right),
 \end{gather}
 where $T_p(\lambda)$ is the time delay operator for coclosed forms defined
 by
 \begin{gather}
  T_p(\lambda)=-\I S_p(\lambda)^{-1} S'_{p-1}(\lambda).
 \end{gather}
 This operator describes the time-delay of coclosed forms of energy 
 $\lambda < \mu_1$ scattered in $M$.  
 In the physics literature the  time delay of the $\ell=0$
 partial wave for potential scattering in $\R^3$ is called the scattering length 
(e.g. \cite{RS79}).
 As the $\ell=0$ partial wave corresponds to the constant function on the sphere at infinity
 we call $T_p(0)$ in analogy with this the scattering length.
 Note however, that in the physics literature the time-delay
 is usually considered for non-relativistic Schr\"odinger mechanics so that
 it differs from the relativistic time-delay by an energy dependent factor. 
 A simple relation that equates time-delay and scattering length can for dimensional reasons
 only hold in relativistic theories. 
\bigskip

\bibliographystyle{amsalpha}

\end{document}